\documentclass[12pt]{amsart}

\usepackage{amsmath}
\usepackage{amscd}
\usepackage{amssymb}
\usepackage{latexsym}
\usepackage{color}

\usepackage[arrow,curve,matrix,tips,2cell]{xy}

\usepackage{multirow}

\newtheorem{theorem}{Theorem}[section]
\newtheorem*{theorem*}{Theorem}
\newtheorem{lemma}[theorem]{Lemma}
\newtheorem*{lemma*}{Lemma}
\newtheorem{corollary}[theorem]{Corollary}
\newtheorem*{corollary*}{Corollary}
\newtheorem{proposition}[theorem]{Proposition}

\newtheorem{remark}[theorem]{Remark}
\newtheorem{question}[theorem]{Question}
\newtheorem{definition}[theorem]{Definition}

\newtheorem{example}[theorem]{Example}
\newcommand{\bgl}{\begin{equation}} 
\newcommand{\egl}{\end{equation}}
\newcommand{\bgloz}{\begin{equation*}} 
\newcommand{\egloz}{\end{equation*}}
\newcommand{\bgln}{\begin{eqnarray}} 
\newcommand{\egln}{\end{eqnarray}}
\newcommand{\bglnoz}{\begin{eqnarray*}} 
\newcommand{\eglnoz}{\end{eqnarray*}}
\newcommand{\btheo}{\begin{theorem}}
\newcommand{\etheo}{\end{theorem}}
\newcommand{\btheooz}{\begin{theorem*}}
\newcommand{\etheooz}{\end{theorem*}}
\newcommand{\blemma}{\begin{lemma}}
\newcommand{\elemma}{\end{lemma}}
\newcommand{\blemmaoz}{\begin{lemma*}}
\newcommand{\elemmaoz}{\end{lemma*}}
\newcommand{\bproof}{\begin{proof}}
\newcommand{\eproof}{\end{proof}}
\newcommand{\bbew}{\begin{beweis}}
\newcommand{\ebew}{\end{beweis}}
\newcommand{\bremark}{\begin{remark}\em}
\newcommand{\eremark}{\end{remark}}
\newcommand{\bquestion}{\begin{question}\em}
\newcommand{\equestion}{\end{question}}
\newcommand{\bdefin}{\begin{definition}}
\newcommand{\edefin}{\end{definition}}
\newcommand{\bprop}{\begin{proposition}}
\newcommand{\eprop}{\end{proposition}}
\newcommand{\bcor}{\begin{corollary}}
\newcommand{\ecor}{\end{corollary}}
\newcommand{\bcoroz}{\begin{corollary*}}
\newcommand{\ecoroz}{\end{corollary*}}
\newcommand{\bfa}{\begin{cases}} 
\newcommand{\efa}{\end{cases}}
\newcommand{\bexample}{\begin{example}\em}
\newcommand{\eexample}{\end{example}}
\newcommand{\cD}{\mathcal D}

\newcommand{\cG}{\mathcal G}

\newcommand{\cL}{\mathcal L}

\newcommand{\cO}{\mathcal O}

\newcommand{\cQ}{\mathcal Q}

\def\Cz{\mathbb{C}}

\def\Nz{\mathbb{N}}

\def\Qz{\mathbb{Q}}
\def\Rz{\mathbb{R}}

\def\Zz{\mathbb{Z}}

\def\1z{\mathbb{1}}

\newcommand{\lori}{\longrightarrow}

\newcommand{\ma}{\mapsto} 
\newcommand\onto{\twoheadrightarrow} 
\newcommand\into{\hookrightarrow} 
\newcommand{\LRarr}{\Leftrightarrow} 

\def\SEMI{\mbox{$\times\kern-2pt\vrule height5pt width.6pt \kern3pt $}}

\newcommand{\Ker}{{\rm Ker\,}}

\newcommand{\Spec}{{\rm Spec\,}} 

\newcommand{\rk}{{\rm rk\,}}

\newcommand{\Coker}{{\rm Coker}\,}
\newcommand{\id}{{\rm id}}

\newcommand{\abs}[1]{\left|#1\right|} 
\newcommand{\norm}[1]{\left\|#1\right\|} 
\newcommand{\defeq}{\mathrel{:=}} 

\newcommand{\dop}{\text{: }} 
\newcommand{\ilim}{\varinjlim} 
\newcommand{\dotcup}{\ensuremath{\mathaccent\cdot\cup}} 

\newcommand{\lge}{\left\{} 
\newcommand{\rge}{\right\}} 
\newcommand{\lru}{\left(} 
\newcommand{\rru}{\right)} 
\newcommand{\lsp}{\left\langle} 
\newcommand{\rsp}{\right\rangle} 
\newcommand{\rukl}[1]{\lru #1 \rru} 
\newcommand{\gekl}[1]{\lge #1 \rge} 
\newcommand{\spkl}[1]{\lsp #1 \rsp} 

\newcommand{\menge}[2]{\gekl{ #1 \dop #2 }} 
\begin{document}

\title[A new approach to constructions of C*-algebras from modular index theory]{A new approach to recent constructions of C*-algebras from modular index theory}

\author{Xin Li}

\address{Xin Li, School of Mathematical Sciences, Queen Mary University of London, Mile End Road, London E1 4NS, UK}
\email{xin.li@qmul.ac.uk}


\begin{abstract}
We present a new approach to C*-algebras recently constructed in the context of modular index theory by Carey, Phillips, Putnam and Rennie. It turns out that their constructions can be identified with full corners of ideals in semigroup C*-algebras. This new point of view leads to a systematic analysis of these algebras and their canonical maximal abelian subalgebras.
\end{abstract}

\maketitle


\setlength{\parindent}{0pt} \setlength{\parskip}{0.5cm}

\section{Introduction}

This paper is about two families of C*-algebras recently constructed in the context of modular index theory \cite{CPPR}. It turns out that these new constructions can be described - in a very concrete and natural way - as certain ideals in semigroup C*-algebras. This leads to a new approach to the C*-algebras from \cite{CPPR} based on recent advances in the study of semigroup C*-algebras \cite{CEL1,CEL2,Li1,Li2}. On the one hand, we are able to explain and simplify structural results from \cite{CPPR}. We provide a better understanding for the results in \cite{CPPR} by relating them to existing ones and by giving conceptual proofs. On the other hand, our approach via semigroup C*-algebras allows us to considerably generalise and strengthen several results from \cite{CPPR}, for instance concerning C*-algebraic structure, or K-theory. Moreover, the C*-algebras constructed in \cite{CPPR} come with canonical maximal abelian subalgebras, which are Cartan subalgebras in the sense of \cite{R08}. We show that there are cases where two C*-algebras (from the families constructed in \cite{CPPR}) are abstractly isomorphic, but not isomorphic as Cartan pairs, i.e., it is impossible to find an isomorphism preserving the canonical Cartan subalgebras.

Let us provide more detailed explanations and formulate our main theorems. First of all, the two classes of C*-algebras constructed in \cite{CPPR} are denoted by $F^{\lambda}$ and $\cQ^{\lambda}$, where the parameter $\lambda$ is a positive real number. The C*-algebras $F^{\lambda}$ and $\cQ^{\lambda}$ are defined as full corners in crossed products, i.e., they are of the form $F^{\lambda} \defeq e(C_0^{\lambda}(\Rz) \rtimes \Gamma_{\lambda})e$ and $\cQ^{\lambda} \defeq e(C_0^{\lambda}(\Rz) \rtimes G_{\lambda})e$. The group $\Gamma_{\lambda}$ is given by $\Gamma_{\lambda} \defeq (\Zz[\lambda,\lambda^{-1}],+)$, i.e., we take the additive group of the subring $\Zz[\lambda,\lambda^{-1}]$ of $\Rz$ generated by $\lambda$ and $\lambda^{-1}$. By construction, $\lambda$ is a unit in $\Zz[\lambda,\lambda^{-1}]$, so we may form the semidirect product $G_{\lambda} \defeq \Gamma_{\lambda} \rtimes \spkl{\lambda}$ with respect to the multiplicative action of $\spkl{\lambda} \defeq \menge{\lambda^i}{i \in \Zz}$ on $\Gamma_{\lambda}$. Moreover, $C_0^{\lambda}(\Rz)$ is the sub-C*-algebra of $\cL^{\infty}(\Rz)$ generated by all characteristic functions on intervals of the form $[a,b)$ for $a, b \in \Gamma_{\lambda}$, $a < b$. The obvious actions of $\Gamma_{\lambda}$ by additive translations and of $G_{\lambda}$ by affine transformations give rise to our crossed products. Finally, $e \in C_0^{\lambda}(\Rz)$ is the characteristic function on $[0,1)$.

The main motivation to introduce and study these C*-algebras $F^{\lambda}$ and $\cQ^{\lambda}$ came from modular index theory. In \cite{CPPR}, an analytic pairing between modular spectral triples and modular $K_1$ is studied in the case of $\cQ^{\lambda}$. We refer the reader to \cite{CPPR} for more details.

Let us form the semigroups $\Gamma_{\lambda}^+ \defeq \Gamma_{\lambda} \cap [0,\infty)$ and $P_{\lambda} \defeq \Gamma_{\lambda}^+ \rtimes \spkl{\lambda}$. The corresponding semigroup C*-algebras $C^*(\Gamma_{\lambda}^+)$ and $C^*(P_{\lambda})$ are generated by isometries representing semigroup elements. By sending these isometries to their unitary counterparts in the corresponding group C*-algebras, we obtain canonical quotient homomorphisms $C^*(\Gamma_{\lambda}^+) \onto C^*(\Gamma_{\lambda})$ and $C^*(P_{\lambda}) \onto C^*(G_{\lambda})$. Let $I(\Gamma_{\lambda}^+)$ and $I(P_{\lambda}))$ be the kernels, i.e., we have short exact sequences $0 \to I(\Gamma_{\lambda}^+) \to C^*(\Gamma_{\lambda}^+) \to C^*(\Gamma_{\lambda}) \to 0$ and $0 \to I(P_{\lambda}) \to C^*(P_{\lambda}) \to C^*(G_{\lambda}) \to 0$. Here is the connection between $F^{\lambda}$, $\cQ^{\lambda}$ and semigroup C*-algebras:

\btheo
\label{THM_FQ-sgpC}
$F^{\lambda}$ is canonically isomorphic to a full corner of $I(\Gamma_{\lambda}^+)$, and $\cQ^{\lambda}$ is canonically isomorphic to a full corner of $I(P_{\lambda})$.
\etheo

We have the following structural results for $F^{\lambda}$ and $\cQ^{\lambda}$:
\btheo
\label{THM_F}
For every $\lambda \in (0,\infty)$, $\lambda \neq 1$, $F^{\lambda}$ is a unital, simple, nuclear, separable C*-algebra with a unique tracial state $\tau$. $F^{\lambda}$ satisfies the UCT. The K-theory of $F^{\lambda}$ is given by
\bgloz
  (K_0(F^{\lambda}),[1]_0,K_1(F^{\lambda})) \cong (K_1(C^*(\Gamma_{\lambda})),[U_1]_1,K_0(C^*(\Gamma_{\lambda}) / \Zz [1]_0).
\egloz
Here $U_1$ is the canonical unitary representing $1 \in \Gamma_{\lambda}$ in $C^*(\Gamma_{\lambda})$. Moreover, the unique trace $\tau$ on $F^{\lambda}$ satisfies $\tau(K_0(F^{\lambda})) = \Gamma_{\lambda}$.

In particular, $F^{\lambda} \cong F^{\mu}$ implies $\Gamma_{\lambda} = \Gamma_{\mu}$.
\etheo
When we write $\Gamma_{\lambda} = \Gamma_{\mu}$, we mean equality as subgroups of $\Rz$.

The structural properties of $F^{\lambda}$ and the K-theory computation also appear in \cite{CPPR}. However, the position of the unit in $K_0$ and the image of the trace have not been determined in \cite{CPPR}. Actually, we will prove this theorem for arbitrary countable dense subgroups of $\Rz$.

\btheo
\label{THM_Q}
For every $\lambda \in (0,\infty)$, $\lambda \neq 1$, $\cQ^{\lambda}$ is a unital UCT Kirchberg algebra. Its K-theory is given by
\bgloz
  (K_0(\cQ^{\lambda}),[1]_0,K_1(\cQ^{\lambda})) \cong (K_1(C^*(\Gamma_{\lambda})) / \Zz [U_{(0,\lambda)}]_1,[U_{(1,1)}]_1^{\bullet},K_0(C^*(\Gamma_{\lambda}) / \Zz [1]_0).
\egloz
Here we write elements of $G_{\lambda}$ in the form $(a,\lambda^i)$, with $a \in \Gamma_{\lambda}$ and $i \in \Zz$. By $[U_{(1,1)}]_1^{\bullet}$, we mean the class of $[U_{(1,1)}]_1$ in the quotient $K_1(C^*(\Gamma_{\lambda})) / \Zz [U_{(0,\lambda)}]_1$.
\etheo

The observation that $\cQ^{\lambda}$ is a UCT Kirchberg algebra also appears in \cite{CPPR}, but our general K-theoretic formula is new.

If $\lambda$ is an algebraic integer, we can make our K-theory computations more precise.
\btheo
\label{THM_Q-algint}
Let $\lambda \in (0,\infty)$, $\lambda \neq 1$ be an algebraic integer. Then
\bgloz
  (K_0(\cQ^{\lambda}),[1]_0,K_1(\cQ^{\lambda})) \cong
  \rukl{\Gamma_{\lambda} / (1 - \lambda) \Gamma_{\lambda} \oplus
  \bigoplus_{j=1}^{\infty} H_{2j+1}(G_{\lambda}), 
  (1,0), 
  \bigoplus_{j=1}^{\infty} H_{2j}(G_{\lambda})}.
\egloz
Here $H_k(G_{\lambda})$ denotes the $k$-th group homology of $G_{\lambda}$. Moreover, $1$ stands for the class of $1 \in \Gamma_{\lambda}$ in $\Gamma_{\lambda} / (1 - \lambda) \Gamma_{\lambda}$.
\etheo
Since group homology can be computed explicitly (as we will explain), this theorem really allows for concrete K-theoretic computations for $\cQ^{\lambda}$. In \cite{CPPR}, such computations (without the position of the unit in $K_0$) have only been carried out for particular values of $\lambda$, under the assumption that both $\lambda$ and $\lambda^{-1}$ are algebraic integers. 

Now let us turn to Cartan pairs. The canonical Cartan subalgebra of $\cQ^{\lambda}$ mentioned at the beginning of the introduction is given by $\cD^{\lambda} \defeq e C_0^{\lambda}(\Rz) e \subseteq e (C_0^{\lambda}(\Rz) \rtimes G_{\lambda}) e = \cQ^{\lambda}$.
\btheo
\label{THM_Q-D}
For all $\lambda \in (0,\infty)$, $\lambda \neq 1$, $\cD^{\lambda}$ is a Cartan subalgebra of $\cQ^{\lambda}$, in the sense of \cite{R08}. Moreover, let $\lambda$ and $\mu$ be two algebraic integers with $\lambda \neq 1 \neq \mu$. If $(\cQ^{\lambda},\cD^{\lambda}) \cong (\cQ^{\mu},\cD^{\mu})$, then $\Gamma_{\lambda} / (1 - \lambda) \Gamma_{\lambda} \cong \Gamma_{\mu} / (1 - \mu) \Gamma_{\mu}$ and $H_i(G_{\lambda}) \cong H_i(G_{\mu})$ for all $i \geq 2$.
\etheo

An isomorphism $(\cQ^{\lambda},\cD^{\lambda}) \cong (\cQ^{\mu},\cD^{\mu})$ means an isomorphism of Cartan pairs, i.e., an isomorphism $\Phi: \: \cQ^{\lambda} \overset{\cong}{\lori} \cQ^{\mu}$ with $\Phi(\cD^{\lambda}) = \cD^{\mu}$. We refer the reader to \cite{R08} for more information about Cartan subalgebras.

We mention that it is possible to formulate and prove similar results for $F^{\lambda}$ as well.

It turns out that we can indeed find algebraic integers $\lambda \neq 1 \neq \mu$ such that
$$
  (K_0(\cQ^{\lambda}),[1]_0,K_1(\cQ^{\lambda})) \cong (K_0(\cQ^{\mu}),[1]_0,K_1(\cQ^{\mu})),
$$
but the condition
$$
  \Gamma_{\lambda} / (1 - \lambda) \Gamma_{\lambda} \cong \Gamma_{\mu} / (1 - \mu) \Gamma_{\mu} \ and \ H_i(\Gamma_{\lambda}) \cong H_i(\Gamma_{\mu}) \ for \ all \ i \geq 2
$$
does not hold. Hence, Theorem~\ref{THM_Q-algint} and Theorem~\ref{THM_Q-D}, together with the Kirchberg-Phillips classification theorem \cite[Chapter~8]{Ror}, imply that $\cQ^{\lambda} \cong \cQ^{\mu}$ but $(\cQ^{\lambda},\cD^{\lambda}) \ncong (\cQ^{\mu},\cD^{\mu})$. For instance, we show that this is the case when the minimal polynomial of $\lambda$ is $T^2 - 3T + 1$ and the minimal polynomial of $\mu$ is $T^3 + T^2 - 1$.

The present paper is structured as follows:
\setlength{\parindent}{15pt} \setlength{\parskip}{0cm}

First, we explain the original construction of $F^{\lambda}$ and $\cQ^{\lambda}$. We also build the bridge to semigroup C*-algebras (Theorem~\ref{THM_FQ-sgpC}), and to groupoid C*-algebras. Furthermore, we prove some first structural properties for $F^{\lambda}$ and $\cQ^{\lambda}$ (Section~\ref{OrgConst-NewDes}).

Secondly, we study $F^{\lambda}$ in more detail and prove Theorem~\ref{THM_F} (Section~\ref{1stF}).

Thirdly, we turn to $\cQ^{\lambda}$ and prove Theorem~\ref{THM_Q}.

Next, we discuss the case where $\lambda$ is an algebraic integer. We prove Theorem~\ref{THM_Q-algint} and Theorem~\ref{THM_Q-D} in Section~\ref{GenCom} and present explicit computations in Section~\ref{ConCom}.

Finally, we close with a list of open questions which arise naturally (Section~\ref{OQ}).
\setlength{\parindent}{0pt} \setlength{\parskip}{0.5cm}

\section{Original constructions and new descriptions}
\label{OrgConst-NewDes}

\subsection{Original constructions}

We briefly recall the construction in \cite{CPPR}. Let $\lambda$ be a positive real number. Let $\Zz[\lambda,\lambda^{-1}]$ be the subring of $\Rz$ generated by $\lambda$ and $\lambda^{-1}$. We denote by $\Gamma_{\lambda}$ the additive group $(\Zz[\lambda,\lambda^{-1}],+)$. The original condition $0 < \lambda < 1$ in \cite{CPPR} is not necessary, at least for the construction of our C*-algebras. But we keep in mind that $\Gamma_{\lambda} = \Gamma_{\lambda^{-1}}$. Let $C_0^{\lambda}(\Rz)$ be the sub-C*-algebra $C^*(\menge{1_{[a,b)}}{a,b \in \Gamma_{\lambda}, a < b})$ of $\cL^{\infty}(\Rz)$. Here $1_{[a,b)}$ is the characteristic function on $[a,b)$. $\Gamma_{\lambda}$ acts on $\cL^{\infty}(\Rz)$ by additive translations, and it is easy to see that $C_0^{\lambda}(\Rz)$ is $\Gamma_{\lambda}$-invariant. Denoting the $\Gamma_{\lambda}$-action on $C_0^{\lambda}(\Rz)$ by $a.f$ for $a \in \Gamma_{\lambda}$ and $f \in C_0^{\lambda}(\Rz)$, we have $(a.f)(x) = f(-a+x)$ for all $x \in \Rz$. Let $A_0^{\lambda} \defeq C_0^{\lambda}(\Rz) \rtimes \Gamma_{\lambda}$ be the corresponding crossed product.

Next, let us form the semidirect product $G_{\lambda} \defeq \Gamma_{\lambda} \rtimes \spkl{\lambda}$ with respect to the multiplicative action of $\spkl{\lambda} \defeq \menge{\lambda^i}{i \in \Zz}$ on $\Gamma_{\lambda}$. Elements of $G_{\lambda}$ are of the form $(a,\lambda^i)$ for $a \in \Gamma_{\lambda}$, $i \in \Zz$, and $(a,\lambda^i)(b,\lambda^j) = (a + \lambda^i b, \lambda^{i+j})$. $G_{\lambda}$ acts on $\cL^{\infty}(\Rz)$ by affine transformations, and $C_0^{\lambda}(\Rz)$ is $G_{\lambda}$-invariant. Denoting the $G_{\lambda}$-action on $C_0^{\lambda}(\Rz)$ by $(a,\lambda^i).f$ for $(a,\lambda^i) \in G_{\lambda}$ and $f \in C_0^{\lambda}(\Rz)$, we have $((a,\lambda^i).f)(x) = f(\lambda^{-i}(-a+x))$. Let $A^{\lambda} \defeq C_0^{\lambda}(\Rz) \rtimes G_{\lambda}$ be the corresponding crossed product.

As the next step, let $e \defeq 1_{[0,1)}$. As $0, 1 \in \Gamma_{\lambda}$, we know that $e \in C_0^{\lambda}(\Rz)$, and we can form the corners $F^{\lambda} \defeq e A_0^{\lambda} e$, $\cQ^{\lambda} \defeq e A^{\lambda} e$. These are the C*-algebras of main interest in this paper. We will focus on their structural properties, and we will also have a look at them from the point of view of Cartan pairs. Apart from that, these C*-algebras are also interesting from the perspective of modular index theory (see \cite{CPPR}).

\subsection{Description as full corners in ideals of semigroup C*-algebras}
\label{corner-ideal-sgpC}

Let us fix a parameter $\lambda \in (0,\infty)$. Set $\Gamma_{\lambda}^+ \defeq \Gamma_{\lambda} \cap [0,\infty)$ and $P_{\lambda} \defeq \Gamma_{\lambda}^+ \rtimes \spkl{\lambda}$. $\Gamma_{\lambda}^+$ and $P_{\lambda}$ are subsemigroups of $\Gamma_{\lambda}$ and $G_{\lambda}$. For simplicity, we drop the index $\lambda$ and write $\Gamma^+$ for $\Gamma_{\lambda}^+$, $\Gamma$ for $\Gamma_{\lambda}$, $P$ for $P_{\lambda}$ and $G$ for $G_{\lambda}$.

The semigroup C*-algebra of $\Gamma^+$ is given by $C^*(\Gamma^+) = C^*(\menge{V_a}{a \in \Gamma^+}) \subseteq \cL(\ell^2 \Gamma^+)$, where $V_a: \: \ell^2 \Gamma^+ \to \ell^2 \Gamma^+, \, \delta_x \ma \delta_{a+x}$. Here $\menge{\delta_x}{x \in \Gamma^+}$ is the canonical orthonormal basis of $\ell^2 \Gamma^+$. Similarly, the semigroup C*-algebra of $P$ is $C^*(P) = C^*(\menge{V_p}{p \in P}) \subseteq \cL(\ell^2 P)$. Note that we write $C^*$ instead of $C^*_{\lambda}$ because in our cases, full and reduced versions coincide. As $\Gamma^+$ is abelian, hence (left) amenable, there exists a homomorphism $C^*(\Gamma^+) \to C^*(\Gamma)$ sending $V_a$ to $U_a$. Here $U_a$ is the canonical unitary in $C^*(\Gamma)$ representing $a$. If we call the kernel of this homomorphism $I(\Gamma^+)$, then we get a short exact sequence
\bgl
\label{ses_ICC-Gamma+}
0 \to I(\Gamma^+) \to C^*(\Gamma^+) \to C^*(\Gamma) \to 0.
\egl
Moreover, let $D(\Gamma^+) \defeq C^*(\menge{1_{a + \Gamma^+}}{a \in \Gamma^+})$ be the canonical commutative sub-C*-algebra of $C^*(\Gamma^+)$. Here $1_{a + \Gamma^+}$ is the characteristic function of $a + \Gamma^+ \subseteq \Gamma^+$, viewed as an element in $\ell^{\infty}(\Gamma^+) \subseteq \cL(\ell^2 \Gamma^+)$. It is easy to see that $D(\Gamma^+) \to \Cz, \, 1_{a + \Gamma^+} \ma 1$ defines a homomorphism. Denoting its kernel by $J(\Gamma^+)$, we obtain a short exact sequence
\bgl
\label{ses_JDC-Gamma+}
  0 \to J(\Gamma^+) \to D(\Gamma^+) \to \Cz \to 0.
\egl
$\Gamma^+$ is abelian, hence left Ore, and its enveloping group is $\Gamma$. So $\Gamma^+ \subseteq \Gamma$ is Toeplitz, in the sense of \cite[Definition~4.1]{Li2}. As in \cite[Definition~3.4]{Li2}, we form $D_{\Gamma^+ \subseteq \Gamma} \defeq C^*(\menge{1_{a + \Gamma^+}}{a \in \Gamma}) \subseteq \cL(\ell^2 \Gamma)$. $\Gamma$ acts on $D_{\Gamma^+ \subseteq \Gamma}$ via $\alpha_a(1_{x + \Gamma^+}) = 1_{a + x + \Gamma^+}$. \cite[Lemma~3.9]{Li2} implies that $C^*(\Gamma^+) = 1_{\Gamma^+} (D_{\Gamma^+ \subseteq \Gamma} \rtimes \Gamma) 1_{\Gamma^+}$. Moreover, the short exact sequences \eqref{ses_JDC-Gamma+} and \eqref{ses_ICC-Gamma+} induce short exact sequences (see \cite[\S~7.1]{Li2})
\bgln
\label{ses_JDC-Gamma}
  && 0 \to J_{\Gamma^+ \subseteq \Gamma} \to D_{\Gamma^+ \subseteq \Gamma} \to \Cz \to 0, \\
\label{ses_JGDGCG-Gamma}
  && 0 \to J_{\Gamma^+ \subseteq \Gamma} \rtimes \Gamma \to D_{\Gamma^+ \subseteq \Gamma} \rtimes \Gamma \to C^*(\Gamma) \to 0.
\egln
It is easy to see that $J_{\Gamma^+ \subseteq \Gamma} = C^*(\menge{1_{a + \Gamma^+} - 1_{b + \Gamma^+}}{a, b \in \Gamma, \, a < b})$ (compare \cite[Definition~7.6]{Li2}). \cite[Proposition~7.9]{Li2} tells us that $I(\Gamma^+) = 1_{\Gamma^+} (J_{\Gamma^+ \subseteq \Gamma} \rtimes \Gamma) 1_{\Gamma^+}$.

In complete analogy, we obtain the following short exact sequences for $P \subseteq G$:
\bgln
\label{ses_ICC-P}
  && 0 \to I(P) \to C^*(P) \to C^*(G) \to 0, \\
\label{ses_JDC-P}
  && 0 \to J(P) \to D(P) \to \Cz \to 0, \\
\label{ses_JDC-G}
  && 0 \to J_{P \subseteq G} \to D_{P \subseteq G} \to \Cz \to 0, \\
\label{ses_JGDGCG-G}
  && 0 \to J_{P \subseteq G} \rtimes G \to D_{P \subseteq G} \rtimes G \to C^*(G) \to 0.
\egln
Again, it is easy to see that $J_{P \subseteq G} = C^*(\menge{1_{(a,1)P} - 1_{(b,1)P}}{a, b \in \Gamma, \, a < b})$. Since $(a,\lambda^i)P = (a,1)P$ for all $i \in \Zz$, it suffices to consider $(a,1)P$. Furthermore, $C^*(P) = 1_P (D_{P \subseteq G} \rtimes G) 1_P$ and $I(P) = 1_P (J_{P \subseteq G} \rtimes G) 1_P$.

To build the bridge between $F^{\lambda}$, $\cQ^{\lambda}$ and semigroup C*-algebras, we need the following
\blemma
\label{totord-proj}
Let $I$ be a totally ordered set, $A$ a C*-algebra generated by projections $e_i$, $i \in I$. Assume that $e_i \lneq e_j$ if $i \lneq j$. Then whenever $B$ is a C*-algebra generated by projections $f_i$, $i \in I$, such that $f_i \leq f_j$ for $i \leq j$, there exists a unique homomorphism $A \to B$ sending $e_i$ to $f_i$ for all $i \in I$.
\elemma
\bproof
This is a special case of \cite[Proposition~2.4]{Li2}. The reason is that $\menge{e_i}{i \in I}$ is obviously a multiplicative semigroup which is independent in the sense of \cite[Definition~2.3]{Li2}.
\eproof

Consider now $C^{\lambda}(\Rz) \defeq C^*({1_{[a,\infty)}}{a \in \Gamma}) \subseteq \cL^{\infty}(\Rz)$. Obviously, $C_0^{\lambda}(\Rz) \subseteq C^{\lambda}(\Rz)$. $\Gamma$ and $G$ act on $C^{\lambda}(\Rz)$ by additive translations and affine transformations, respectively.
\bcor
There exists a $\Gamma$-equivariant isomorphism $D_{\Gamma^+ \subseteq \Gamma} \overset{\cong}{\lori} C^{\lambda}(\Rz)$ sending $1_{a + \Gamma^+}$ to $1_{[a,\infty)}$ for all $a \in \Gamma$. Also, there exists a $G$-equivariant isomorphism $D_{P \subseteq G} \overset{\cong}{\lori} C^{\lambda}(\Rz)$ sending $1_{(a,1)P}$ to $1_{[a,\infty)}$ for all $a \in \Gamma$.
\ecor
\bproof
This is a direct consequence of Lemma~\ref{totord-proj}.
\eproof

\bcor
\label{J-C0l}
There is a $\Gamma$-equivariant isomorphism $J_{\Gamma^+ \subseteq \Gamma} \overset{\cong}{\lori} C_0^{\lambda}(\Rz)$ sending $1_{a + \Gamma^+} - 1_{b + \Gamma^+}$ to $1_{[a,b)}$ for all $a, b \in \Gamma$, $a < b$. Moreover, $1_{\Gamma^+} - 1_{1 + \Gamma^+}$ is a full projection in $I(\Gamma^+)$ and in $J_{\Gamma^+ \subseteq \Gamma}$, $J_{\Gamma^+ \subseteq \Gamma} \rtimes \Gamma \cong C_0^{\lambda}(\Rz) \rtimes \Gamma$, and $(1_{\Gamma^+} - 1_{1 + \Gamma^+})I(\Gamma^+)(1_{\Gamma^+} - 1_{1 + \Gamma^+}) = (1_{\Gamma^+} - 1_{1 + \Gamma^+})(J_{\Gamma^+ \subseteq \Gamma} \rtimes \Gamma)(1_{\Gamma^+} - 1_{1 + \Gamma^+}) \cong F^{\lambda}$.

There is also a $G$-equivariant isomorphism $J_{P \subseteq G} \overset{\cong}{\lori} C_0^{\lambda}(\Rz)$ sending $1_{(a,1)P} - 1_{(b,1)P}$ to $1_{[a,b)}$ for all $a, b \in \Gamma$, $a < b$. Moreover, $1_P - 1_{(1,1)P}$ is a full projection in $I(P)$ and in $J_{P \subseteq G}$, $J_{P \subseteq G} \rtimes G \cong C_0^{\lambda}(\Rz) \rtimes G$, and $(1_P - 1_{(1,1)P})I(P)(1_P - 1_{(1,1)P}) = (1_P - 1_{(1,1)P})(J_{P \subseteq G} \rtimes G)(1_P - 1_{(1,1)P}) \cong \cQ^{\lambda}$.
\ecor
\bproof
Just note that the isomorphisms $D_{\Gamma^+ \subseteq \Gamma} \cong C^{\lambda}(\Rz)$ and $D_{P \subseteq G} \cong C^{\lambda}(\Rz)$ send $J_{\Gamma^+ \subseteq \Gamma}$ to $C_0^{\lambda}(\Rz)$, $1_{\Gamma^+} - 1_{1 + \Gamma^+}$ to $e = 1_{[0,1)}$, $J_{P \subseteq G}$ to $C_0^{\lambda}(\Rz)$, and $1_P - 1_{(1,1)P}$ to $e$. Moreover, it is easy to see that $1_{\Gamma^+} - 1_{1 + \Gamma^+}$ and $1_P - 1_{(1,1)P}$ are full projections in the corresponding C*-algebras (see also Lemma~\ref{N-full} and Remark~\ref{Rem_N-full}).
\eproof

This proves Theorem~\ref{THM_FQ-sgpC}.

\bremark
$\Gamma$ could be any subgroup of $(\Rz,+)$, not necessarily of the form $\Gamma_{\lambda}$, and all our constructions make sense, and analogous results as above hold. In particular, $1_{\Gamma^+} - 1_{a + \Gamma^+}$ is a full projection in the corresponding C*-algebras whenever $a \in \Gamma$ satisfies $a > 0$.
\eremark

\subsection{Description as groupoid C*-algebras}
\label{des-gpd}

We describe all the C*-algebras from the previous section, in particular $F^{\lambda}$ and $\cQ^{\lambda}$, as groupoid C*-algebras. These descriptions follow from general considerations in \cite[\S~5.4, \S~5.5, \S~7.1]{Li2}, and we refer the reader to \cite{Li2} for details.

Let us consider the slightly more general situation of an arbitrary subgroup $\Gamma$ of $(\Rz,+)$. Set $\Gamma^+ \defeq \Gamma \cap [0,\infty)$. Without loss of generality, we may assume that $\Gamma$ is dense in $\Rz$, since otherwise $(\Gamma,\Gamma^+)$ is isomorphic to $(\Zz,\Nz)$, and this case is well-understood. Let $\Omega(\Gamma^+) \defeq \Spec(D(\Gamma^+))$ and $\Omega_{\Gamma^+ \subseteq \Gamma} \defeq \Spec(D_{\Gamma^+ \subseteq \Gamma})$. We have $\Omega(\Gamma^+) \subseteq \Omega_{\Gamma^+ \subseteq \Gamma}$. Moreover, let $\chi_{\infty}$ be the character on $D(\Gamma^+)$ satisfying $\chi_{\infty}(1_{a + \Gamma^+}) = 1$ for all $a \in \Gamma^+$. Set $O(\Gamma^+) \defeq \Spec(D(\Gamma^+)) \setminus \gekl{\chi_{\infty}}$ and $O_{\Gamma^+ \subseteq \Gamma} \setminus \gekl{\chi_{\infty}}$. The canonical identifications $D(\Gamma^+) \cong C(\Omega(\Gamma^+))$ and $D_{\Gamma^+ \subseteq \Gamma} \cong C_0(\Omega_{\Gamma^+ \subseteq \Gamma})$ restrict to isomorphisms $J(\Gamma^+) \cong C_0(O(\Gamma^+))$ and $J_{\Gamma^+ \subseteq \Gamma} \cong C_0(O_{\Gamma^+ \subseteq \Gamma})$. These isomorphisms induce isomorphisms
\bglnoz
  && D_{\Gamma^+ \subseteq \Gamma} \rtimes \Gamma 
  \cong C^*(\Gamma \ltimes \Omega_{\Gamma^+ \subseteq \Gamma}), 
  \qquad
  C^*(\Gamma^+) \cong C^*(\Gamma \ltimes \Omega_{\Gamma^+ \subseteq \Gamma} \big \vert {}_{\Omega(\Gamma^+)}^{\Omega(\Gamma^+)}), \\
  && J_{\Gamma^+ \subseteq \Gamma} \rtimes \Gamma 
  \cong C^*(\Gamma \ltimes O_{\Gamma^+ \subseteq \Gamma}),
  \qquad I(\Gamma^+) \cong C^*(\Gamma \ltimes O_{\Gamma^+ \subseteq \Gamma} \big \vert {}_{O(\Gamma^+)}^{O(\Gamma^+)}),
\eglnoz
where $\Gamma \ltimes \Omega_{\Gamma^+ \subseteq \Gamma}$ and $\Gamma \ltimes O_{\Gamma^+ \subseteq \Gamma}$ are the transformation groupoids attached to $\Gamma \curvearrowright \Omega_{\Gamma^+ \subseteq \Gamma}$ and $\Gamma \curvearrowright O_{\Gamma^+ \subseteq \Gamma}$, the dual actions of $\Gamma \curvearrowright D_{\Gamma^+ \subseteq \Gamma}$ and $\Gamma \curvearrowright J_{\Gamma^+ \subseteq \Gamma}$. Moreover, it is easy to check that $\Omega(\Gamma^+)$ is a $\Gamma \ltimes \Omega_{\Gamma^+ \subseteq \Gamma}$-full compact open subset of $\Omega_{\Gamma^+ \subseteq \Gamma}$, meaning that $\Omega(\Gamma^+)$ meets every $\Gamma$-orbit in $\Omega_{\Gamma^+ \subseteq \Gamma}$. Similarly, $O(\Gamma^+)$ is a $\Gamma \ltimes O_{\Gamma^+ \subseteq \Gamma}$-full clopen subset of $O_{\Gamma^+ \subseteq \Gamma}$. Now choose $a \in \Gamma$ with $a > 0$. Then $N(\Gamma^+) \defeq \menge{\chi \in \Omega(\Gamma^+)}{\chi(1_{a + \Gamma^+}) = 0}$ is a compact open subset of $\Omega(\Gamma^+)$ with $1_{N(\Gamma^+)} = 1_{\Gamma^+} - 1_{a + \Gamma^+}$.

\blemma
\label{N-full}
$N(\Gamma^+)$ is $\Gamma \ltimes O_{\Gamma^+ \subseteq \Gamma}$-full in $O_{\Gamma^+ \subseteq \Gamma}$.
\elemma
\bproof
Take $\chi \in O_{\Gamma^+ \subseteq \Gamma}$. Since $\chi \neq \chi_{\infty}$, there exists $a \in \Gamma$ such that $\chi(1_{a + \Gamma^+}) = 0$. There also exists $b \in \Gamma$ such that $\chi(1_{b + \Gamma^+}) = 1$. Moreover, for all $b, c \in \Gamma$ with $b \geq c$, $\chi(1_{b + \Gamma^+}) = 1$ implies $\chi(1_{c + \Gamma^+}) = 1$. Therefore, $\menge{a \in \Gamma}{\chi(1_{a + \Gamma^+}) = 0}$ is bounded below (by $b$). Let $\tau \defeq \inf \menge{a \in \Gamma}{\chi(1_{a + \Gamma^+}) = 0}$. Since $\Gamma$ is dense in $\Rz$, there exists $b \in \Gamma$ with $\tau \in (b,a+b]$. Hence $((-b).\chi)(1_{\Gamma^+}) = \chi(1_{b + \Gamma^+}) = 1$ and $((-b).\chi)(1_{a + \Gamma^+}) = \chi(1_{a + b + \Gamma^+}) = 0$. Therefore, $(-b).\chi \in N(\Gamma^+)$.
\eproof

\bremark
\label{Rem_N-full}
Lemma~\ref{N-full} implies that $1_{\Gamma^+} - 1_{a + \Gamma^+}$ is a full projection in 
$J_{\Gamma^+ \subseteq \Gamma} \rtimes \Gamma$.
\eremark

Obviously, with $E \defeq 1_{\Gamma^+} - 1_{a + \Gamma^+}$, we have $E(J_{\Gamma^+ \subseteq \Gamma} \rtimes \Gamma)E \cong C^*(\Gamma \ltimes O_{\Gamma^+ \subseteq \Gamma} \big \vert {}_{N(\Gamma^+)}^{N(\Gamma^+)})$. In particular, for $\Gamma = \Gamma_{\lambda}$ and $a = 1$, we obtain $N(\Gamma_{\lambda}^+) = \menge{\chi \in \Omega(\Gamma_{\lambda}^+)}{\chi(1_{1 + \Gamma_{\lambda}^+}) = 0}$ and
$$
  F^{\lambda} 
  \cong
  C^*(\Gamma_{\lambda} \ltimes O_{\Gamma_{\lambda}^+ \subseteq \Gamma_{\lambda}} 
  \big \vert {}_{N(\Gamma_{\lambda}^+)}^{N(\Gamma_{\lambda}^+)}).$$

Completely analogous descriptions exist in the case of $P_{\lambda} \subseteq G_{\lambda}$ instead of $\Gamma^+ \subseteq \Gamma$. As in the previous section, we fix $\lambda \in (0,\infty)$ and write $P$ for $P_{\lambda}$, $G$ for $G_{\lambda}$. Let $\Omega(P) \defeq \Spec(D(P))$ and $\Omega_{P \subseteq G} \defeq \Spec(D_{P \subseteq G})$. Denote by $\chi_{\infty}$ the character on $D(P)$ satisfying $\chi_{\infty}(1_{pP}) = 1$ for all $p \in P$. Set $O(P) \defeq \Omega(P) \setminus \gekl{\chi_{\infty}}$ and $O_{P \subseteq G} \defeq \Omega_{P \subseteq G} \setminus \gekl{\chi_{\infty}}$. Also, let $N(P) \defeq \menge{\chi \in \Omega(P)}{\chi(1_{(1,1)P}) = 1}$. As before, $\Omega(P)$ is a $G \ltimes \Omega_{P \subseteq G}$-full compact open subset of $\Omega_{P \subseteq G}$, $O(P)$ is a $G \ltimes O_{P \subseteq G}$-full clopen subset of $O_{P \subseteq G}$, and $N(P)$ is a $G \ltimes O_{P \subseteq G}$-full compact open subset of $O_{P \subseteq G}$. Furthermore, we have canonical identifications
\bglnoz
  && D_{P \subseteq G} \rtimes G \cong C^*(G \ltimes \Omega_{P \subseteq G}),
  \qquad
  C^*(P) \cong C^*(G \ltimes \Omega_{P \subseteq G} \big \vert {}_{\Omega(P)}^{\Omega(P)}), \\
  && J_{P \subseteq G} \rtimes G
  \cong C^*(G \ltimes O_{P \subseteq G}),
  \qquad
  I(P) \cong C^*(G \ltimes O_{P \subseteq G} \big \vert {}_{O(P)}^{O(P)}), \\
  &{\rm and}& \cQ^{\lambda} \cong C^*(G \ltimes O_{P \subseteq G} \big \vert {}_{N(P)}^{N(P)}).
\eglnoz

\subsection{Structural properties}

With the help of the groupoid description from the previous section, it is easy to derive structural properties for $E(J_{\Gamma^+ \subseteq \Gamma} \rtimes \Gamma)E$, in particular $F^{\lambda}$, and $\cQ^{\lambda}$.

\blemma
\label{nuc-UCT}
For every dense subgroup $\Gamma$ of $\Rz$, $E(J_{\Gamma^+ \subseteq \Gamma} \rtimes \Gamma)E$ is nuclear and satisfies the UCT. Hence for every $\lambda \in (0,\infty)$, $\lambda \neq 1$, $F^{\lambda}$ is nuclear and satisfies the UCT.

For all $\lambda \in (0,\infty)$, $\lambda \neq 1$, $\cQ^{\lambda}$ is nuclear and satisfies the UCT.
\elemma
\bproof
We have seen that both $E(J_{\Gamma^+ \subseteq \Gamma} \rtimes \Gamma)E$ and $\cQ^{\lambda}$ are isomorphic to groupoid C*-algebras. The underlying groupoids are amenable since $\Gamma$ and $G_{\lambda}$ are amenable groups. Our lemma follows from \cite{AR} and \cite{Tu}.
\eproof

\blemma
For every dense subgroup $\Gamma$ of $\Rz$, the actions $\Gamma \curvearrowright \Omega_{\Gamma^+ \subseteq \Gamma}$ and $\Gamma \curvearrowright O_{\Gamma^+ \subseteq \Gamma}$ are topologically free.
\elemma
\bproof
Since $O_{\Gamma^+ \subseteq \Gamma}$ is dense in $\Omega_{\Gamma^+ \subseteq \Gamma}$, it suffices to show that $\Gamma \curvearrowright O_{\Gamma^+ \subseteq \Gamma}$ is topologically free. First of all, it follows from \cite[\S~2.3]{Li2} that basic open sets in $O_{\Gamma^+ \subseteq \Gamma}$ are of the form $U_{a,b} = \menge{\chi \in \Omega_{\Gamma^+ \subseteq \Gamma}}{\chi(1_{a + \Gamma^+}) = 1, \, \chi(1_{b + \Gamma^+}) = 0}$ for $a,b \in \Gamma$, $a < b$. For $x \in \Gamma$, let $\chi_x \in O_{\Gamma^+ \subseteq \Gamma}$ be the character determined by
$$
  \chi_x(1_{y + \Gamma^+}) = 
  \bfa
    1 & if \ x + \Gamma^+ \subseteq y + \Gamma^+ \LRarr x \geq y, \\
    0 & else.
  \efa
$$
Obviously, $\chi_a \in U_{a,b}$, so that $\menge{\chi_a}{a \in \Gamma}$ is dense in $O_{\Gamma^+ \subseteq \Gamma}$. But $\Gamma$ acts freely on $\menge{\chi_a}{a \in \Gamma}$. Therefore, $\Gamma \curvearrowright O_{\Gamma^+ \subseteq \Gamma}$ is topologically free.
\eproof

\blemma
\label{G-action-TF}
For every $\lambda \in (0,\infty)$, $\lambda \neq 1$, $G_{\lambda} \curvearrowright \Omega_{P_{\lambda} \subseteq G_{\lambda}}$ and $G_{\lambda} \curvearrowright O_{P_{\lambda} \subseteq G_{\lambda}}$ are topologically free.
\elemma
\bproof
As before, we fix $\lambda \in (0,\infty)$, $\lambda \neq 1$, and write $\Gamma$ for $\Gamma_{\lambda}$, $P$ for $P_{\lambda}$ and $G$ for $G_{\lambda}$. As for the previous lemma, it suffices to show that $G_{\lambda} \curvearrowright O_{P_{\lambda} \subseteq G_{\lambda}}$ is topologically free. By \cite[\S~2.3]{Li2}, basic open sets in $O_{P \subseteq G}$ are given by $U_{a,b} = \menge{\chi \in \Omega_{P \subseteq G}}{\chi(1_{(a,1)P}) = 1, \, \chi(1_{(b,1)P}) = 0}$ for $a,b \in \Gamma$, $a < b$. Again, for $x \in \Gamma$, let $\chi_x \in O_{P \subseteq G}$ be the character determined by
$$
  \chi_x(1_{(y,1)P}) = 
  \bfa
    1 & if \ (x,1)P \subseteq (y,1)P \LRarr x \geq y, \\
    0 & else.
  \efa
$$
Now let $(c,\lambda^i) \in G$. If $(c,\lambda^i).\chi = \chi$, then $c + \lambda^i x = x$, and this implies $(1 - \lambda^i)x = c$. If $(c,\lambda^i)$ is not the identity element, i.e., $(c,\lambda^i) \neq (0,1)$, then we see that $(c,\lambda^i)$ fixes at most one character in $\menge{\chi_x}{x \in \Gamma}$. At the same time, for all $a, b \in \Gamma$ with $a < b$, there exist infinitely many $x \in \Gamma$ such that $\chi_x \in U_{a,b}$. Namely, since $\lambda \neq 1$, $\Gamma$ is dense in $\Rz$, so that there exist infinitely many $x \in \Gamma$ with $a < x < b$, and for every such $x$, we have $\chi_x \in U_{a,b}$. All in all, we deduce that for all $(0,1) \neq (c,\lambda^i) \in G$ and for all $a, b \in \Gamma$ with $a < b$, there exists $x \in \Gamma$ with $\chi_x \in U_{a,b}$ and $(c,\lambda^i).\chi \neq \chi$.
\eproof

\bcor
\label{gpd-TF}
For every dense subgroup $\Gamma$ of $\Rz$, the groupoids
$$
  \Gamma \ltimes \Omega_{\Gamma^+ \subseteq \Gamma} \big \vert {}_{\Omega(\Gamma^+)}^{\Omega(\Gamma^+)}, \, \Gamma \ltimes O_{\Gamma^+ \subseteq \Gamma} \big \vert {}_{O(\Gamma^+)}^{O(\Gamma^+)} \ and \ \Gamma \ltimes O_{\Gamma^+ \subseteq \Gamma} \big \vert {}_{N(\Gamma^+)}^{N(\Gamma^+)}
$$
are topologically principal, in the sense of \cite{R08}.

For every $\lambda \in (0,\infty)$, $\lambda \neq 1$, the groupoids
$$
  G_{\lambda} \ltimes \Omega_{P_{\lambda} \subseteq G_{\lambda}} \big \vert {}_{\Omega(P_{\lambda})}^{\Omega(P_{\lambda})}, \ G_{\lambda} \ltimes O_{P_{\lambda} \subseteq G_{\lambda}} \big \vert {}_{O(P_{\lambda})}^{O(P_{\lambda})} \ and \ G_{\lambda} \ltimes O_{P_{\lambda} \subseteq G_{\lambda}} \big \vert {}_{N(P_{\lambda})}^{N(P_{\lambda})}
$$
are topologically principal.
\ecor
\bproof
This follows from the general (and obvious) observation that restrictions of topologically principal groupoids to open subsets of their unit spaces are again topologically principal.
\eproof

Let $\cD^{\lambda}$ be the commutative sub-C*-algebra $e C_0^{\lambda}(\Rz) e$ of $\cQ^{\lambda}$.
\bcor
For every dense subgroup $\Gamma$ of $\Rz$, $(C^*(\Gamma^+),D(\Gamma^+))$, $(I(\Gamma^+),J(\Gamma^+))$ and $(E(J_{\Gamma^+ \subseteq \Gamma} \rtimes \Gamma)E,E(J_{\Gamma^+ \subseteq \Gamma})E)$ are Cartan pairs (in the sense of \cite[Definition~5.1]{R08}). Therefore, for all $\lambda \in (0,\infty)$, $\lambda \neq 1$, $(F^{\lambda},eC_0^{\lambda}(\Rz)e)$ is a Cartan pair.

For all $\lambda \in (0,\infty)$, $\lambda \neq 1$, $(C^*(P_{\lambda}),D(P_{\lambda}))$, $(I(P_{\lambda}),J(P_{\lambda}))$ and
$$
  (\cQ^{\lambda},\cD^{\lambda}) \cong ((1_{P_{\lambda}} - 1_{(1,1)P_{\lambda}})(J_{P_{\lambda} \subseteq G_{\lambda}} \rtimes \Gamma)(1_{P_{\lambda}} - 1_{(1,1)P_{\lambda}}),(1_{P_{\lambda}} - 1_{(1,1)P_{\lambda}})(J_{P_{\lambda} \subseteq G_{\lambda}})(1_{P_{\lambda}} - 1_{(1,1)P_{\lambda}}))
$$
are Cartan pairs.
\ecor
\bproof
This follows from Corollary~\ref{gpd-TF} and \cite[Theorem~5.2]{R08}.
\eproof

\blemma
Let $\Gamma$ be a dense subgroup of $\Rz$. If $I$ is a nonzero ideal of $C^*(\Gamma^+)$, then there exists $0 \neq a \in \Gamma^+$ with $1_{\Gamma^+} - 1_{a + \Gamma^+} \in I$.

Let $\lambda \in (0,\infty)$, $\lambda \neq 1$. If $I$ is a nonzero ideal of $C^*(P_{\lambda})$, then there exists $0 \neq a \in \Gamma^+$ with $1_{P_{\lambda}} - 1_{(a,1)P_{\lambda}} \in I$.
\elemma
\bproof
We give the proof for $C^*(\Gamma^+)$, the one for $C^*(P_{\lambda})$ is completely analogous. First of all, since $\Gamma \curvearrowright \Omega_{\Gamma^+ \subseteq \Gamma}$ is topologically free, \cite{AS} tells us that every nonzero ideal of $D_{\Gamma^+ \subseteq \Gamma} \rtimes \Gamma$ has nonzero intersection with $D_{\Gamma^+ \subseteq \Gamma} \cong C_0(\Omega_{\Gamma^+ \subseteq \Gamma})$. Since $C^*(\Gamma^+)$ is (isomorphic to) the full corner $1_{\Gamma^+}(D_{\Gamma^+ \subseteq \Gamma} \rtimes \Gamma)1_{\Gamma^+}$, and we know that $D(\Gamma^+) = 1_{\Gamma^+}D_{\Gamma^+ \subseteq \Gamma} 1_{\Gamma^+}$, $C^*(\Gamma^+)$ also has the property that every nonzero ideal intersects $D(\Gamma^+)$ non-trivially. But by \cite[\S~2.6]{Li1}, every nonzero ideal of $D(\Gamma^+)$ must contain an element of the form $1_{\Gamma^+} - 1_{a + \Gamma^+} \in I$ for some $0 \neq a \in \Gamma^+$.
\eproof

\blemma
For every dense subgroup $\Gamma$ of $\Rz$, $I(\Gamma^+)$ is simple. Also, for all $\lambda \in (0,\infty)$, $\lambda \neq 1$, $I(P_{\lambda})$ is simple.
\elemma
\bproof
Again, we treat the case of $I(\Gamma^+)$, the proof for $I(P_{\lambda})$ is completely analogous. Let $(0) \neq I \triangleleft I(\Gamma^+)$ be a nonzero ideal of $I(\Gamma^+)$. As $I$ is also a nonzero ideal of $C^*(\Gamma^+)$, the previous lemma tells us that there exists $0 \neq a \in \Gamma^+$ with $1_{\Gamma^+} - 1_{a + \Gamma^+} \in I$. Let $b \in \Gamma^+$ be arbitrary, and choose $n \in \Nz$ such that $b \leq na$. Then $1_{b + \Gamma^+} \geq 1_{na + \Gamma^+}$, so we conclude that
\bgloz
  1_{\Gamma^+} - 1_{b + \Gamma^+} \leq 1_{\Gamma^+} - 1_{na + \Gamma^+} = \sum_{m=0}^{n-1} 1_{ma + \Gamma^+} - 1_{(m+1)a + \Gamma^+} = \sum_{m=0}^{n-1} V_{ma} (1_{\Gamma^+} - 1_{a + \Gamma^+}) V_{ma}^* \in I.
\egloz
Thus $1_{\Gamma^+} - 1_{b + \Gamma^+} \in I$. But $b \in \Gamma^+$ was arbitrary, and $\menge{1_{\Gamma^+} - 1_{b + \Gamma^+}}{b \in \Gamma^+}$ generates $I(\Gamma^+)$. Hence we must have $I = I(\Gamma^+)$.
\eproof

\bremark
For $\Gamma^+$, this observation appears in \cite{Mur}.
\eremark

Since simplicity is preserved under Morita equivalence, we obtain
\bcor
\label{simple}
For every dense subgroup $\Gamma$ of $\Rz$, $J_{\Gamma^+ \subseteq \Gamma} \rtimes \Gamma$ and $E(J_{\Gamma^+ \subseteq \Gamma} \rtimes \Gamma)E$ are simple. Hence for every $\lambda \in (0,\infty)$, $\lambda \neq 1$, $F^{\lambda}$ is simple.

Moreover, for every $\lambda \in (0,\infty)$, $\lambda \neq 1$, $J_{P_{\lambda} \subseteq G_{\lambda}} \rtimes G_{\lambda}$ and
$\cQ^{\lambda} \cong (1_{P_{\lambda}} - 1_{(1,1)P_{\lambda}})(J_{P_{\lambda} \subseteq G_{\lambda}} \rtimes \Gamma)(1_{P_{\lambda}} - 1_{(1,1)P_{\lambda}})$ are simple.
\ecor

\section{First family of C*-algebras}
\label{1stF}

In this section, let $\Gamma$ be a dense subgroup of $\Rz$. Moreover, take $a \in \Gamma$ with $a > 0$ and write $E$ for $1_{\Gamma^+} - 1_{a + \Gamma^+}$.

\bprop
\label{JGamma-tau}
$E(J_{\Gamma^+ \subseteq \Gamma} \rtimes \Gamma)E$ has a unique tracial state $\tau$ given by
\bgl
\label{tau}
  \tau((1_{b + \Gamma^+} - 1_{c + \Gamma^+})U_{\gamma})
  = 
  \bfa
  0 & if \ \gamma \neq 0, \\
  \tfrac{1}{a}(c-b) & if \ \gamma = 0,
  \efa
\egl
for all $b, c \in \Gamma^+$, $b < c$ and for all $\gamma \in \Gamma$.

In particular, for every $\lambda \in (0,\infty)$, $\lambda \neq 1$, $F^{\lambda}$ has a unique tracial state.
\eprop
\bproof
We first construct a tracial state with the desired properties. First of all, the C*-algebra $C^{\lambda}(\Rz)$ (see \S~\ref{corner-ideal-sgpC}) is obviously contained in the C*-algebra of bounded Borel functions on $\Rz$. Under the isomorphism $J_{\Gamma^+ \subseteq \Gamma} \cong C^{\lambda}(\Rz)$ in Corollary~\ref{J-C0l}, $E = 1_{\Gamma^+} - 1_{a + \Gamma^+}$ is mapped to $1_{[0,a)}$, so that $E J_{\Gamma^+ \subseteq \Gamma} E$ is identified with a sub-C*-algebra of bounded Borel functions on $[0,a)$. Hence Lebesgue measure, scaled by $\tfrac{1}{a}$, gives rise to a (tracial) state $\mu$ on $E J_{\Gamma^+ \subseteq \Gamma} E$. Furthermore, let $\theta: \: J_{\Gamma^+ \subseteq \Gamma} \rtimes \Gamma \to J_{\Gamma^+ \subseteq \Gamma}$ be the canonical faithful conditional expectation. Set $\tau \defeq \mu \circ (\theta \vert_{E(J_{\Gamma^+ \subseteq \Gamma} \rtimes \Gamma)E})$. It is clear that $\tau$ is a state, and $\tau$ obviously satisfies \eqref{tau}. An easy computation using \eqref{tau} shows that $\tau$ has the trace property.

We now prove uniqueness. Let $\tau'$ be another tracial state on $E(J_{\Gamma^+ \subseteq \Gamma} \rtimes \Gamma)E$. Let $\gamma \in \Gamma$ be nonzero. For $b, c \in \Gamma^+$, $b < c$, let $b_i$ ($i = 0, \dotsc, n+1$) be elements of $\Gamma^+$ with $b = b_0 < b_1 < \dotso < b_n < b_{n+1} = c$ such that $b_{i+1} - b_i < \abs{\gamma}$. (Here we use that $\Gamma$ is dense in $\Rz$.) Then $\tau'((1_{b + \Gamma^+} - 1_{c + \Gamma^+})U_{\gamma}) = \sum_{i=0}^n \tau'((1_{b_i + \Gamma^+} - 1_{b_{i+1} + \Gamma^+})U_{\gamma})$. Now $[b_i,b_{i+1}) \cap [\gamma + b_i,\gamma + b_{i+1}) = \emptyset$ since $b_{i+1} - b_i < \abs{\gamma}$. Hence
\bglnoz
  \tau'((1_{b_i + \Gamma^+} - 1_{b_{i+1} + \Gamma^+})U_{\gamma})
  &=& \tau'(
  (1_{b_i + \Gamma^+} - 1_{b_{i+1} + \Gamma^+})
  U_{\gamma}
  (1_{b_i + \Gamma^+} - 1_{b_{i+1} + \Gamma^+})
  ) \\
  &=& \tau'(
  (1_{b_i + \Gamma^+} - 1_{b_{i+1} + \Gamma^+})
  (1_{\gamma + b_i + \Gamma^+} - 1_{\gamma + b_{i+1} + \Gamma^+})
  U_{\gamma}
  )
  = \tau'(0) = 0.  
\eglnoz
Now let $\mu'$ be the composition $1_{[0,a)} C_0^{\lambda}(\Rz) 1_{[0,a)} \cong E J_{\Gamma^+ \subseteq \Gamma} E \overset{\tau'}{\lori} \Cz$. The trace property of $\tau'$ implies that $\mu'(1_{[b,c)}) = \mu'(1_{[\gamma + b,\gamma + c)})$ for all $b, c \in \Gamma^+$, $b < c \leq a$ and for all $\gamma \in \Gamma$ such that $0 \leq \gamma + b$, $\gamma + c \leq a$. As $\Gamma$ is dense in $\Rz$, this shows that $\mu'$ is induced by a translation invariant measure on $[0,a)$. Therefore, $\mu'$ must be given by Lebesgue measure, scaled by $\tfrac{1}{a}$.
\eproof

We turn to K-theory.
\bprop
\label{K-Gamma}
For every countable subgroup $\Gamma$ of $\Rz$, the boundary maps in the six term exact sequence for $0 \to J_{\Gamma^+ \subseteq \Gamma} \rtimes \Gamma \to D_{\Gamma^+ \subseteq \Gamma} \rtimes \Gamma \to C^*(\Gamma) \to 0$ (see \eqref{ses_JGDGCG-Gamma}) induce isomorphisms $K_1(C^*(\Gamma)) \overset{\cong}{\lori} K_0(J_{\Gamma^+ \subseteq \Gamma} \rtimes \Gamma)$ and $K_0(C^*(\Gamma)) / \Zz [1]_0 \overset{\cong}{\lori} K_1(J_{\Gamma^+ \subseteq \Gamma} \rtimes \Gamma)$. The first isomorphism sends $[1_{\Gamma^+} - 1_{a + \Gamma^+}]_0$ to $[U_a]_1$. Moreover, $\Zz [1]_0$ is a direct summand (isomorphic to $\Zz$) in $K_0(C^*(\Gamma))$.
\eprop
\bproof
It is easy to see that \cite[Corollary~3.14]{CEL2} applies and yields in our situation that the homomorphism $\Cz \to D_{\Gamma^+ \subseteq \Gamma} \rtimes \Gamma, \, z \ma z \cdot 1_{\Gamma^+}$ induces an isomorphism $K_*(\Cz) \cong K_*(D_{\Gamma^+ \subseteq \Gamma} \rtimes \Gamma)$. Moreover, $\Cz \to C^*(\Gamma), \, z \ma z \cdot 1$ induces a split-injective homomorphism in K-theory because a split is provided by $C^*(\Gamma) \to \Cz, \, U_{\gamma} \ma 1$. Analysing the six term exact sequence for $0 \to J_{\Gamma^+ \subseteq \Gamma} \rtimes \Gamma \to D_{\Gamma^+ \subseteq \Gamma} \rtimes \Gamma \to C^*(\Gamma) \to 0$ with the help of these observations, we arrive at the desired conclusions.
\eproof

As $I(\Gamma^+)$ and $E(J_{\Gamma^+ \subseteq \Gamma} \rtimes \Gamma)E$ are full corners in $J_{\Gamma^+ \subseteq \Gamma} \rtimes \Gamma$, we obtain
\bcor
For every countable, dense subgroup $\Gamma$ of $\Rz$,
$$
  K_i(I(\Gamma^+)) \cong K_i(E(J_{\Gamma^+ \subseteq \Gamma} \rtimes \Gamma)E)
  \cong K_i(J_{\Gamma^+ \subseteq \Gamma} \rtimes \Gamma) \cong
  \bfa
    K_1(C^*(\Gamma)) & if \ i = 0, \\
    K_0(C^*(\Gamma)) / \Zz [1]_0 & if \ i = 1.
  \efa
$$
In particular, we have for all $\lambda \in (0,\infty)$, $\lambda \neq 1$,
$$
  (K_0(F^{\lambda}),[1]_0,K_1(F^{\lambda})) 
  \cong (K_1(C^*(\Gamma_{\lambda})),[U_1]_1,K_0(C^*(\Gamma_{\lambda})) / \Zz [1]_0).
$$
\ecor 

\bremark
The last result has been proven (apart from the statement concerning the class of the unit) in \cite{CPPR} with different methods.
\eremark

Let $\tau$ be the unique tracial state on $E(J_{\Gamma^+ \subseteq \Gamma} \rtimes \Gamma)E$, where $E = 1_{\Gamma^+} - 1_{a + \Gamma^+}$. Our next goal is to compute $\tau(K_0(E(J_{\Gamma^+ \subseteq \Gamma} \rtimes \Gamma)E))$. Let $\Gamma$ be a countable, dense subgroup of $\Rz$. Then the representation
\bgl
\label{rep_R-Gamma}
  C_0(\Rz \cup \gekl{+ \infty}) \to \cL(\ell^2 \Gamma), \, f \ma M_f,
\egl
where $M_f \delta_x = f(x) \delta_x$, is faithful. Approximating (in $\norm{\cdot}_{\infty}$) every $f \in C_0(\Rz \cup \gekl{+ \infty})$ by linear combinations of $1_{[a,\infty)}$, $a \in \Gamma$, we see that the representation \eqref{rep_R-Gamma} gives rise to a $\Gamma$-equivariant embedding $C_0(\Rz \cup \gekl{+ \infty}) \into D_{\Gamma^+ \subseteq \Gamma}$. This embedding restricts to a $\Gamma$-equivariant embedding $C_0(\Rz) \into J_{\Gamma^+ \subseteq \Gamma}$.

\bprop
\label{R-JGamma}
For every countable, dense subgroup $\Gamma$ of $\Rz$, $C_0(\Rz) \rtimes \Gamma \into J_{\Gamma^+ \subseteq \Gamma} \rtimes \Gamma$ induces a $K_0$-isomorphism and an isomorphism $K_1(C_0(\Rz) \rtimes \Gamma) / K_1(C_0(\Rz)) \overset{\cong}{\lori} K_1(J_{\Gamma^+ \subseteq \Gamma} \rtimes \Gamma)$.
\eprop
\bproof
Consider the following commutative diagram with exact rows:
\bgl
\label{CD-exrows}
  \xymatrix{
  0 \ar[r] & C_0(\Rz) \rtimes \Gamma \ar[d] \ar[r] &
  C_0(\Rz \cup \gekl{+ \infty}) \rtimes \Gamma \ar[d] \ar[r] &
  C^*(\Gamma) \ar[d]^{=} \ar[r] & 0   
  \\
  0 \ar[r] & J_{\Gamma^+ \subseteq \Gamma} \rtimes \Gamma \ar[r] &
  D_{\Gamma^+ \subseteq \Gamma} \rtimes \Gamma \ar[r] & 
  C^*(\Gamma) \ar[r] & 0
  }
\egl
Since $C_0(\Rz \cup \gekl{+ \infty}) \sim_h \gekl{0}$ and $\Gamma$ is a torsion-free group satisfying the Baum-Connes conjecture, principle (DC) from \cite[\S~1]{CEL2} implies that $K_*(C_0(\Rz \cup \gekl{+ \infty}) \rtimes \Gamma) \cong \gekl{0}$. Hence \eqref{CD-exrows} gives in K-theory commuting squares
\bgl
\label{CD1}
  \xymatrix{
  K_1(C^*(\Gamma)) \ar[d]^{=} \ar[r]^-{\cong} & K_0(C_0(\Rz) \rtimes \Gamma) \ar[d]
  \\
  K_1(C^*(\Gamma)) \ar[r]^-{\cong} & K_0(J_{\Gamma^+ \subseteq \Gamma} \rtimes \Gamma)
  }
\egl
\bgl
\label{CD2}
  \xymatrix{
  K_0(C^*(\Gamma)) \ar[d]^{=} \ar[r]^-{\cong} & K_1(C_0(\Rz) \rtimes \Gamma) \ar[d]
  \\
  K_0(C^*(\Gamma)) \ar[r] & K_1(J_{\Gamma^+ \subseteq \Gamma} \rtimes \Gamma)
  }
\egl
In \eqref{CD1}, the bottom arrow is an isomorphism by Proposition~\ref{K-Gamma}. Hence $C_0(\Rz) \rtimes \Gamma \into J_{\Gamma^+ \subseteq \Gamma} \rtimes \Gamma$ induces a $K_0$-isomorphism. By Proposition~\ref{K-Gamma}, we know that the kernel of the bottom arrow in \eqref{CD2} is $\Zz [1]_0$. So it remains to determine the image of $\Zz [1]_0$ under $K_0(C^*(\Gamma)) \overset{\cong}{\lori} K_1(C_0(\Rz) \rtimes \Gamma)$. Consider the following commutative diagram with exact rows:
\bgloz
  \xymatrix{
  0 \ar[r] & C_0(\Rz) \ar@{^{(}->}[d] \ar[r] &
  C_0(\Rz \cup \gekl{+ \infty}) \ar@{^{(}->}[d] \ar[r] & \Cz \ar@{^{(}->}[d] \ar[r] & 0   
  \\
  0 \ar[r] & C_0(\Rz) \rtimes \Gamma \ar[r] &
  C_0(\Rz \cup \gekl{+ \infty}) \rtimes \Gamma \ar[r] & C^*(\Gamma) \ar[r] & 0   
  }
\egloz
The vertical arrows are the canonical inclusions. In K-theory, we obtain a commuting square
\bgloz
  \xymatrix{
  K_0(\Cz) \ar[d] \ar[r]^-{\cong} & K_1(C_0(\Rz)) \ar[d]
  \\
  K_0(C^*(\Gamma)) \ar[r]^-{\cong} & K_1(C_0(\Rz) \rtimes \Gamma)
  }
\egloz
This shows that the image of $\Zz [1]_0$ under $K_0(C^*(\Gamma)) \overset{\cong}{\lori} K_1(C_0(\Rz) \rtimes \Gamma)$ coincides with the image of $K_1(C_0(\Rz)) \to K_1(C_0(\Rz) \rtimes \Gamma)$.
\eproof

\bprop
For every countable, dense subgroup $\Gamma$ of $\Rz$, $\tau(K_0(E(J_{\Gamma^+ \subseteq \Gamma} \rtimes \Gamma)E)) = \tfrac{1}{a} \Gamma$.
\eprop
\bproof
By \cite{CZ}, there exists a unique densely defined trace $\tau_{\Gamma^+ \subseteq \Gamma}$ on $J_{\Gamma^+ \subseteq \Gamma} \rtimes \Gamma$ such that $\tau_{\Gamma^+ \subseteq \Gamma} \vert_{E(J_{\Gamma^+ \subseteq \Gamma} \rtimes \Gamma)E} = \tau$. The composition $C_0(\Rz) \rtimes \Gamma \into J_{\Gamma^+ \subseteq \Gamma} \rtimes \Gamma \overset{\tau_{\Gamma^+ \subseteq \Gamma}}{\lori} \Cz$ must be the unique densely defined trace on $C_0(\Rz) \rtimes \Gamma$ induced by Lebesgue measure, scaled by $\tfrac{1}{a}$. Hence by \cite[Lemma~2.2]{JX}, we have $\tau_{\Gamma^+ \subseteq \Gamma}(K_0(C_0(\Rz) \rtimes \Gamma)) = \tfrac{1}{a} \Gamma$. But since $C_0(\Rz) \rtimes \Gamma \into J_{\Gamma^+ \subseteq \Gamma} \rtimes \Gamma$ induces a $K_0$-isomorphism, we deduce that $\tau_{\Gamma^+ \subseteq \Gamma}(K_0(J_{\Gamma^+ \subseteq \Gamma} \rtimes \Gamma)) = \tfrac{1}{a} \Gamma$. Furthermore, the canonical inclusion $E(J_{\Gamma^+ \subseteq \Gamma} \rtimes \Gamma)E \into J_{\Gamma^+ \subseteq \Gamma} \rtimes \Gamma$ induces a $K_0$-isomorphism, and $\tau_{\Gamma^+ \subseteq \Gamma} \vert_{E(J_{\Gamma^+ \subseteq \Gamma} \rtimes \Gamma)E} = \tau$. Thus $\tau(K_0(E(J_{\Gamma^+ \subseteq \Gamma} \rtimes \Gamma)E)) = \tau_{\Gamma^+ \subseteq \Gamma}(K_0(J_{\Gamma^+ \subseteq \Gamma} \rtimes \Gamma)) = \tfrac{1}{a} \Gamma$.
\eproof

For $\Gamma = \Gamma_{\lambda}$ and $a = 1 \in \Gamma_{\lambda}$, we obtain
\bcor
\label{tFGamma}
For every $\lambda \in (0,\infty)$, $\lambda \neq 1$, the unique tracial state $\tau$ on $F^{\lambda}$ satisfies $\tau(F^{\lambda}) = \Gamma_{\lambda}$.
\ecor

\bremark
Proposition~\ref{JGamma-tau} and Proposition~\ref{R-JGamma} have been proven in \cite{JX} by different methods. Moreover, the special case of Proposition~\ref{JGamma-tau} concerning $F^{\lambda}$ has been observed in \cite[\S~2.6]{CPPR}.
\eremark

Note that Theorem~\ref{THM_F} follows from Lemma~\ref{nuc-UCT}, Lemma~\ref{simple}, Proposition~\ref{JGamma-tau}, Proposition~\ref{K-Gamma} and Corollary~\ref{tFGamma}.

\section{Second family of C*-algebras}
\label{2ndF}

In this section, we fix $\lambda \in (0,\infty)$, $\lambda \neq 1$, and set $\Gamma \defeq \Gamma_{\lambda}$, $P \defeq P_{\lambda}$ and $G \defeq G_{\lambda}$.

\bprop
\label{loc-bnd}
$G \curvearrowright O_{P \subseteq G}$ is a topologically free local boundary action (in the sense of \cite{LS}).
\eprop
\bproof
We have already seen in Lemma~\ref{G-action-TF} that $G \curvearrowright O_{P \subseteq G}$ is topologically free. To prove that $G \curvearrowright O_{P \subseteq G}$ is a local boundary action, we have to show (by definition) that for every non-empty open subset $U$ of $O_{P \subseteq G}$, there exists an open subset $\Delta$ of $U$ and $g \in G$ such that $g \overline{\Delta} \subsetneq \Delta$. Basic open sets of $O_{P \subseteq G}$ are of the form $U_{w,z} = \menge{\chi \in \Omega_{P \subseteq G}}{\chi(1_{(w,1)P}) = 1, \, \chi(1_{(z,1)P}) = 0}$ for $w,z \in \Gamma$, $w < z$. Thus we may assume that $U = U_{w,z}$ for some $w, z$. It suffices to find $x, y \in \Gamma$ with $x < y$ and $(a,\lambda^i) \in G$ such that $U_{x,y} \subseteq U_{w,z}$ and $(a,\lambda^i).U_{x,y} \subsetneq U_{x,y}$. (Note that $U_{x,y}$ is clopen.) As $\Gamma$ is dense in $\Rz$, we can always find $x, y \in \Gamma$ with $w < x < y < z$. Choose $i \in \Zz$ such that $0 < \lambda^i < 1$. Using density of $\Gamma$ again, we can find $a \in \Gamma$ such that $(1 - \lambda^i)x < a < (1 - \lambda^i)y$. Obviously, $U_{x,y} \subseteq U_{w,z}$ as $w < x < y < z$. Moreover, we claim that $(a,\lambda^i).U_{x,y} = U_{a + \lambda^i x, a + \lambda^i y}$. Indeed, given $\chi \in \Omega_{P \subseteq G}$, $((a,\lambda^i).\chi)(1_{(v,1)P}) = \chi(1_{(\lambda^{-i}(-a+v),1)P})$. So $((a,\lambda^i).\chi)(1_{(a + \lambda^i v,1)P}) = 1$ if and only if $\chi(1_{(v,1)P}) = 1$. As $x < a + \lambda^i x < a + \lambda^i y < y$ by construction, we conclude that $(a,\lambda^i).U_{x,y} \subsetneq U_{x,y}$.
\eproof

\bcor
\label{Q-Kiralg}
For all $\lambda \in (0,\infty)$, $\lambda \neq 1$, $\cQ^{\lambda}$ is a unital UCT Kirchberg algebra.
\ecor
\bproof
It follows from Lemma~\ref{nuc-UCT} and Lemma~\ref{simple} that $\cQ^{\lambda}$ is simple, separable, unital, nuclear and satisfies the UCT. \cite[Theorem~9]{LS} and Proposition~\ref{loc-bnd} imply that $J_{P \subseteq G} \rtimes G \cong C_0(O_{P \subseteq G}) \rtimes G$ is purely infinite simple. As $\cQ^{\lambda}$ is isomorphic to a full corner of $J_{P \subseteq G} \rtimes G$, it must be purely infinite simple as well.
\eproof

We now compute K-theory for $\cQ^{\lambda}$.
\bprop
The boundary maps in the six term exact sequence for $0 \to J_{P \subseteq G} \rtimes G \to D_{P \subseteq G} \rtimes G \to C^*(G) \to 0$ (see \eqref{ses_JGDGCG-G}) induce isomorphisms $K_1(C^*(G)) / \Zz [U_{(0,\lambda)}]_1 \overset{\cong}{\lori} K_0(J_{P \subseteq G} \rtimes G)$ and $K_0(C^*(G)) / \Zz [1]_0 \overset{\cong}{\lori} K_1(J_{P \subseteq G} \rtimes G)$. Under the first isomorphism, $[1_P - 1_{(1,1)P}]_0$ corresponds to the class of $[U_{(1,1)}]_1 \in K_1(C^*(G))$ in the quotient $K_1(C^*(G)) / \Zz [U_{(0,\lambda)}]_1$. Moreover, $\Zz [U_{(0,\lambda)}]_1$ and $\Zz [1]_0$ are direct summands (isomorphic to $\Zz$) in $K_1(C^*(\Gamma))$ and $K_0(C^*(\Gamma))$, respectively.
\eprop
\bproof
It is easy to see that \cite[Corollary~3.14]{CEL2} applies and yields in our situation that $C^*(\spkl{\lambda}) \to D_{P \subseteq G} \rtimes G, \, U_{\lambda^i} \ma U_{\lambda^i} 1_P$ induces an isomorphism $K_*(C^*(\spkl{\lambda})) \cong K_*(D_{P \subseteq G} \rtimes G)$. Moreover, $\spkl{\lambda} \to G, \, \lambda^i \ma (0,\lambda^i) $ and $G \to \spkl{\lambda}, \, (a,\lambda^i) \ma \lambda^i$ induce homomorphisms $C^*(\spkl{\lambda}) \to C^*(G)$ and $C^*(G) \to C^*(\spkl{\lambda})$ such that the composition $C^*(\spkl{\lambda}) \to C^*(G) \to C^*(\spkl{\lambda})$ is the identity. Therefore, $C^*(\spkl{\lambda}) \to C^*(G)$ induces a split-injective homomorphism in K-theory. With the help of these observations, it is straightforward to analyse the six term exact sequence for $0 \to J_{P \subseteq G} \rtimes G \to D_{P \subseteq G} \rtimes G \to C^*(G) \to 0$ and to derive the desired conclusions.
\eproof

\bremark
Of course, with a different sign convention for the index maps, $[1_P - 1_{(1,1)P}]_0$ would correspond to the class of $- [U_{(1,1)}]_1$.
\eremark

\bcor
\label{K-Q-general}
For every $\lambda \in (0,\infty)$, $\lambda \neq 1$, there is an isomorphism
$$
  (K_0(\cQ^{\lambda}),[1]_0,K_1(\cQ^{\lambda})) 
  \cong (K_1(C^*(G_{\lambda})) / \Zz [U_{(0,\lambda)}]_1,
  [U_{(1,1)}]_1^{\bullet},K_0(C^*(G_{\lambda})) / \Zz [1]_0).
$$
\ecor
Note that by the Kirchberg-Phillips classification theorem \cite[Chapter~8]{Ror}, this is a complete isomorphism invariant for $\cQ^{\lambda}$. Also note that Corollary~\ref{Q-Kiralg} and Corollary~\ref{K-Q-general} prove Theorem~\ref{THM_Q}.

Let us turn to isomorphism of Cartan pairs. Recall that we write $\cD^{\lambda}$ for the sub-C*-algebra $e C_0^{\lambda}(\Rz) e$ of $\cQ^{\lambda}$.
\blemma
\label{criterion_Cartan}
For $\lambda, \mu \in (0,\infty)$ with $\lambda \neq 1 \neq \mu$, $(\cQ^{\lambda},\cD^{\lambda}) \cong (\cQ^{\mu},\cD^{\mu})$ implies that
$$
  H_*(G_{\lambda}, C_c(O_{P_{\lambda} \subseteq G_{\lambda}},\Zz))
  \cong H_*(G_{\lambda}, C_c(O_{P_{\mu} \subseteq G_{\mu}},\Zz)).
$$
\elemma
Here $H_*(G_{\lambda}, C_c(O_{P_{\lambda} \subseteq G_{\lambda}},\Zz))$ is the group homology of $G_{\lambda}$ with coefficients in $C_c(O_{P_{\lambda} \subseteq G_{\lambda}},\Zz)$, where $G_{\lambda} \curvearrowright C_c(O_{P_{\lambda} \subseteq G_{\lambda}},\Zz)$ is induced by $G_{\lambda} \curvearrowright O_{P_{\lambda} \subseteq G_{\lambda}}$.
\bproof
We know by \S~\ref{des-gpd} that
$
  (\cQ^{\lambda},\cD^{\lambda}) \cong 
  \rukl{C^*(G_{\lambda} \ltimes O_{P_{\lambda} \subseteq G_{\lambda}}
  \big \vert {}_{N(P_{\lambda})}^{N(P_{\lambda})}),C(N(P_{\lambda}))}
$, and similarly for $(\cQ^{\mu},\cD^{\mu})$. Hence, by \cite[Proposition~4.13]{R08}, we have $(\cQ^{\lambda},\cD^{\lambda}) \cong (\cQ^{\mu},\cD^{\mu})$ if and only if $G_{\lambda} \ltimes O_{P_{\lambda} \subseteq G_{\lambda}} \big \vert {}_{N(P_{\lambda})}^{N(P_{\lambda})} \cong G_{\mu} \ltimes O_{P_{\mu} \subseteq G_{\mu}} \big \vert {}_{N(P_{\mu})}^{N(P_{\mu})}$ as topological groupoids. The latter statement implies that the groupoid homologies (in the sense of \cite[\S~3.1]{Matui}) must be isomorphic, i.e., $H_*(G_{\lambda} \ltimes O_{P_{\lambda} \subseteq G_{\lambda}} \big \vert {}_{N(P_{\lambda})}^{N(P_{\lambda})},\Zz) \cong H_*(G_{\mu} \ltimes O_{P_{\mu} \subseteq G_{\mu}} \big \vert {}_{N(P_{\mu})}^{N(P_{\mu})},\Zz)$. But \cite[Theorem~3.6~(2)]{Matui} implies that
$
  H_*(G_{\lambda} \ltimes O_{P_{\lambda} \subseteq G_{\lambda}} \big \vert {}_{N(P_{\lambda})}^{N(P_{\lambda})},\Zz) \cong H_*(G_{\lambda} \ltimes O_{P_{\lambda} \subseteq G_{\lambda}},\Zz) \cong H_*(G_{\lambda}, C_c(O_{P_{\lambda} \subseteq G_{\lambda}},\Zz))
$,
and similarly, $H_*(G_{\mu} \ltimes O_{P_{\mu} \subseteq G_{\mu}} \big \vert {}_{N(P_{\mu})}^{N(P_{\mu})},\Zz) \cong H_*(G_{\mu}, C_c(O_{P_{\mu} \subseteq G_{\mu}},\Zz))$.
\eproof

Note that a much more obvious invariant of $(\cQ^{\lambda},\cD^{\lambda})$ (apart from the isomorphism class of $\cQ^{\lambda}$) is the isomorphism class of $\cD^{\lambda}$. However, all the $\cD^{\lambda}$ are mutually isomorphic for $1 \neq \lambda \in (0,\infty)$ since $\cD^{\lambda} \cong C(N(P_{\lambda}))$, and $N(P_{\lambda})$ is homeomorphic to the Cantor space for all $1 \neq \lambda \in (0,\infty)$.

\section{Kirchberg algebras parametrised by algebraic integers}

Our goal is to compute $(K_0(\cQ^{\lambda}),[1]_0,K_1(\cQ^{\lambda}))$ more precisely for algebraic integers $\lambda$. At the same time, we compute $H_*(G_{\lambda}, C_c(O_{P_{\lambda} \subseteq G_{\lambda}},\Zz))$ so that it is easy to compare with K-theory. 

\subsection{General computations}
\label{GenCom}

We need the following
\blemma
\label{GGG}
Let $(G,+)$ be an abelian group, $\varphi: \: G \to G$ an endomorphism. Let $G_{\infty} \defeq \ilim \gekl{G,\varphi}$ be the inductive limit of $G \overset{\varphi}{\lori} G \overset{\varphi}{\lori} G \overset{\varphi}{\lori} \dotso$. Let $\varphi_{\infty}: \: G_{\infty} \to G_{\infty}$ be the map induced by $\varphi$, given by the following morphism of inductive systems
\bgloz
  \xymatrix{
  G \ar[d]^{\varphi} \ar[r]^{\varphi} & G \ar[d]^{\varphi} \ar[r]^{\varphi} &
  G \ar[d]^{\varphi} \ar[r]^{\varphi} & \dotso \\
  G \ar[r]^{\varphi} & G \ar[r]^{\varphi} & G \ar[r]^{\varphi} & \dotso 
  }
\egloz
Let $\iota: \: G \to G_{\infty}$ be the homomorphism of the first copy of $G$ to $G_{\infty}$ (which is part of the structure of the inductive limit). Then $\iota$ induces an isomorphism $\Ker(\id - \varphi) \overset{\cong}{\lori} \Ker(\id - \varphi_{\infty})$ and $\Coker(\id - \varphi) \overset{\cong}{\lori} \Coker(\id - \varphi_{\infty})$.
\elemma
\bproof
We work with the following concrete model for $G_{\infty}$: We write $[g_n]_n$ for the class of $(g_n)_n \in \prod_{n=1}^{\infty} G$ in the quotient $\prod_{n=1}^{\infty} G / \bigoplus_{n=1}^{\infty} G$. Let $G_{\infty}$ be the subgroup of the quotient group $\prod_{n=1}^{\infty} G / \bigoplus_{n=1}^{\infty} G$ generated by all elements $[g_n]_n$ for which there exists $M \in \Nz$ such that $g_{m+1} = \varphi(g_m)$ for all $m \geq M$. Together with the homomorphisms $\iota_n: \: G \to G_{\infty}, \, g \ma [\underbrace{0, \dotsc, 0}_{n-1}, g, \varphi(g), \varphi^2(g), \dotsc]$, $G_{\infty}$ is the inductive limit of $G \overset{\varphi}{\lori} G \overset{\varphi}{\lori} G \overset{\varphi}{\lori} \dotso$. Moreover, $\varphi_{\infty} [g_1,g_2,g_3,\dotsc] = [\varphi(g_1),\varphi(g_2),\varphi(g_3),\dotsc] = [g_2,g_3,g_4,\dotsc]$.

Let us first show that $\iota \defeq \iota_1$ induces an isomorphism $\Ker(\id - \varphi) \overset{\cong}{\lori} \Ker(\id - \varphi_{\infty})$. Take $g \in \Ker(\id - \varphi)$. Then $g = \varphi(g)$, so $\iota(g) = [g,g,g,\dotsc] \in \Ker(\id - \varphi_{\infty})$. If $\iota(g) = 0$, then $(g,g,g,\dotsc)$ lies in $\bigoplus_{n=1}^{\infty} G$, and hence $g$ must be $0$. Therefore, $\iota \vert_{Ker(\id - \varphi)}$ is injective. Now take $[g_n]_n \in \Ker(\id - \varphi_{\infty})$. Then $[g_1,g_2,\dotsc] = \varphi_{\infty} [g_1,g_2,\dotsc] = [g_2,g_3,\dotsc]$. Therefore, there exists $M \in \Nz$ such that $g_m = g_{m+1}$ and $g_{m+1} = \varphi(g_m)$ for all $m \geq M$. This means that $g \defeq g_M$ lies in $\Ker(\id - \varphi)$, and $[g_n]_n = [g,g,\dotsc] = \iota(g)$. Thus $\iota \vert_{Ker(\id - \varphi)}$ is surjective.

Secondly, $\iota$ induces a homomorphism $\iota^{\bullet}: \: \Coker(\id - \varphi) \lori \Coker(\id - \varphi_{\infty})$. The reason is that if $g = h - \varphi(h)$ for some $h \in G$, then $\iota(g) = [g,\varphi(g),\dotsc] = [h - \varphi(h),\varphi(h) - \varphi^2(h),\dotsc] = [h,\varphi(h),\dotsc] - [\varphi(h),\varphi^2(h),\dotsc] = (\id - \varphi_{\infty})(\iota(h))$. To see that $\iota^{\bullet}$ is injective, assume that $\iota(g) = [g,\varphi(g),\dotsc]$ is of the form $(\id - \varphi_{\infty}) [g_1,g_2,\dotsc] = [g_1,g_2,\dotsc] - [g_2,g_3,\dotsc] = [g_1 - g_2,g_2 - g_3,\dotsc]$. Then there exists $M \in \Nz$ such that $\varphi^m(g) = g_m - g_{m+1}$ and $g_{m+1} = \varphi(g_m)$ for all $m \geq M$. Then $g = g - \varphi^M(g) + \varphi^M(g) = (\id - \varphi)(\sum_{m=0}^{M-1} \varphi^m(g)) + g_M - \varphi(g_M) \in (\id - \varphi)(G)$. Let us show that $\iota^{\bullet}$ is surjective. For $[g_n]_n \in G_{\infty}$, there exists $M \in \Nz$ with $g_{m+1} = \varphi(g_m)$ for all $m \geq M$. Then we have modulo $(\id - \varphi_{\infty})(G_{\infty})$:
$
  [g_n]_n = [g_n]_n - (\id - \varphi_{\infty})(\sum_{m=0}^{M-1} \varphi_{\infty}^m [g_n]_n) = [g_n]_n - (\id - \varphi_{\infty}^M) [g_n]_n = \varphi_{\infty}^M [g_n]_n = [g_M,\varphi(g_M),\varphi^2(g_M),\dotsc] = \iota(g_M)
$.
\eproof

We now apply this lemma. Let $\lambda$ be an algebraic integer with minimal polynomial of degree $d$. Let $\varphi: \: \Zz[\lambda] \to \Zz[\lambda]$ be the homomorphism given by multiplication with $\lambda$. Here $\Zz[\lambda]$ is the subring of $\Rz$ generated by $1$ and $\lambda$. Moreover, let $\Lambda^k \varphi: \: \Lambda^k \Zz[\lambda] \to \Lambda^k \Zz[\lambda]$ be the induced homomorphism on the $k$-th exterior power.
\blemma
\label{K-G}
There is an isomorphism $K_0(C^*(G_{\lambda})) \cong \bigoplus_{j=0}^{\infty} \Coker(\id - \Lambda^{2j} \varphi) \oplus \bigoplus_{j=0}^{\infty} \Ker(\id - \Lambda^{2j+1} \varphi)$ such that $\Zz [1]_0$ is identified with $\Coker(\id - \Lambda^0 \varphi) = \Lambda^0 \Zz[\lambda] \cong \Zz$.

We also have an isomorphism $K_1(C^*(G_{\lambda})) \cong \bigoplus_{j=0}^{\infty} \Coker(\id - \Lambda^{2j+1} \varphi) \oplus \bigoplus_{j=0}^{\infty} \Ker(\id - \Lambda^{2j} \varphi)$ sending $\Zz [U_{(0,\lambda)}]_1$ to $\Ker(\id - \Lambda^0 \varphi) = \Lambda^0 \Zz[\lambda] \cong \Zz$ and $[U_{(1,1)}]_1$ to $[1] = [\lambda^0] \in \Coker(\id - \Lambda^1 \varphi) = \Zz[\lambda] / (1 - \lambda) \Zz[\lambda] \cong \Gamma_{\lambda} / (1 - \lambda) \Gamma_{\lambda}$.
\elemma
\bproof
Obviously, we can identify $\Gamma_{\lambda} = \Zz[\lambda,\lambda^{-1}]$ with $\ilim \gekl{\Zz[\lambda],\varphi}$ such that $\varphi_{\infty}$ corresponds to the automorphism of $\Gamma_{\lambda}$ given by multiplication with $\lambda$. Therefore, we obtain $K_*(C^*(\Gamma_{\lambda})) \cong \ilim \gekl{K_*(C^*(\Zz[\lambda])),\varphi_*}$, and $(\varphi_{\infty})_*$ corresponds to $(\varphi_*)_{\infty}$. The Pimsner-Voiculescu exact sequence for $C^*(G_{\lambda}) \cong (C^*(\Gamma_{\lambda}) \rtimes \spkl{\lambda}$ gives rise to short exact sequences
\bglnoz
  && 0 \to \Coker(\id - K_0(\varphi_{\infty})) \to K_0(C^*(G_{\lambda})) \to \Ker(\id - K_1(\varphi_{\infty})) \to 0, \\
  && 0 \to \Coker(\id - K_1(\varphi_{\infty})) \to K_1(C^*(G_{\lambda})) \to \Ker(\id - K_0(\varphi_{\infty})) \to 0.
\eglnoz
Lemma~\ref{GGG} implies that $\Coker(\id - K_i(\varphi_{\infty})) \cong \Coker(\id - K_i(\varphi))$ and $\Ker(\id - K_i(\varphi_{\infty})) \cong \Ker(\id - K_i(\varphi))$ for all $i=0,1$. We know that $\Zz[\lambda] = \bigoplus_{i=0}^{d-1} \Zz \lambda^i \cong \Zz^d$, so that $K_i(C^*(\Zz[\lambda])) \cong K_i(C^*(\Zz^d)) \cong \Zz^{2^{d-1}}$ is a free abelian group. Thus $\Ker(\id - K_i(\varphi))$ is free abelian as well, and 
$$
  K_i(C^*(G_{\lambda})) \cong \Coker(\id - K_i(\varphi)) \oplus \Ker(\id - K_{i+1}(\varphi))
$$
for all $i=0,1$. Moreover, we have natural isomorphisms $K_i(C^*(\Zz[\lambda])) \cong \bigoplus_{j=0}^{\infty} \Lambda^{2j+i} \Zz[\lambda]$ sending $\Zz [1]_0$ to $\Lambda^0 \Zz[\lambda]$ and $[U_1]_1$ to $1 \in \Lambda^1 \Zz[\lambda]$. Naturality means that under these isomorphisms, $K_i(\varphi)$ corresponds to $\bigoplus_{j=0}^{\infty} \Lambda^{2j+i} \varphi$. Finally, an easy computation shows that the index map in the Pimsner-Voiculescu exact sequence for $C^*(G_{\lambda}) \cong C^*(\Gamma_{\lambda}) \rtimes \spkl{\lambda}$ sends $\Zz [U_{(0,\lambda)}]_1 \in K_1(C^*(G_{\lambda}))$ to $\Zz [1]_0 \in K_0(C^*(\Gamma_{\lambda}))$.
\eproof

\bcor
\label{K-Q-algint}
For every algebraic integer $\lambda \in (0,\infty)$, $\lambda \neq 1$,
\bglnoz
  (K_0(\cQ^{\lambda}),[1]_0,K_1(\cQ^{\lambda})) &\cong& 
  \left(
  \bigoplus_{j=0}^{\infty} \Coker(\id - \Lambda^{2j+1} \varphi) \oplus
  \bigoplus_{j=1}^{\infty} \Ker(\id - \Lambda^{2j} \varphi), \right. \\
  && [1] \in \Coker(\id - \Lambda^1 \varphi) 
  \cong \Gamma_{\lambda} / (1 - \lambda) \Gamma_{\lambda}, \\
  && \left. \bigoplus_{j=1}^{\infty} \Coker(\id - \Lambda^{2j} \varphi) \oplus
  \bigoplus_{j=0}^{\infty} \Ker(\id - \Lambda^{2j+1} \varphi)
  \right).
\eglnoz
\ecor
\bproof
This follows from Corollary~\ref{K-Q-general} and Lemma~\ref{K-G}.
\eproof

\bprop
\label{gphom'}
If $\lambda \in (0,\infty)$ is an algebraic integer, then $H_{k+1}(G_{\lambda}) \cong \Coker(\id - \Lambda^{k+1} \varphi) \oplus \Ker(\id - \Lambda^k \varphi)$ for all $k = 0,1,2, \dotsc$.
\eprop
\bproof
The Hochschild-Serre spectral sequence \cite[Chapter~VII, Theorem~(6.3)]{Brown} gives rise to short exact sequences
\bgl
\label{ses-gphom'}
  0 \to \Coker(\id - H_{k+1}(\varphi_{\infty})) \to H_{k+1}(G_{\lambda}) \to \Ker(\id - H_k(\varphi_{\infty})) \to 0.
\egl
As group homology is continuous (see \cite[Theorem~5.5]{Beyl}), an application of Lemma~\ref{GGG} as in our K-theoretic computation yields $\Coker(\id - H_{k+1}(\varphi_{\infty})) \cong \Coker(\id - H_{k+1}(\varphi)) \cong \Coker(\id - \Lambda^{k+1} \varphi)$ and $\Ker(\id - H_k(\varphi_{\infty})) \cong \Ker(\id - H_k(\varphi)) \cong \Ker(\id - \Lambda^k \varphi)$. As $\Lambda^k \Zz[\lambda]$ is free abelian, so is $\Ker(\id - \Lambda^k \varphi)$. Thus \eqref{ses-gphom'} splits. Our proposition follows.
\eproof
Note that Corollary~\ref{K-Q-algint} and Proposition~\ref{gphom'} imply Theorem~\ref{THM_Q-algint}.

Motivated by Lemma~\ref{criterion_Cartan}, we also compute $H_*(G_{\lambda}, C_c(O_{P_{\lambda} \subseteq G_{\lambda}},\Zz))$.
\bprop
\label{gphom}
For all $\lambda \in (0,\infty)$,
\bglnoz
  H_0(G_{\lambda}, C_c(O_{P_{\lambda} \subseteq G_{\lambda}},\Zz))
  &\cong& \Gamma_{\lambda} / (1 - \lambda) \Gamma_{\lambda} \\
  H_k(G_{\lambda}, C_c(O_{P_{\lambda} \subseteq G_{\lambda}},\Zz))
  &\cong& H_{k+1}(G_{\lambda}) \ \ \ \ \ \ for \ all \ k \geq 1.
\eglnoz
\eprop
\bproof
We fix $\lambda \in (0,\infty)$ and write $\Gamma$ for $\Gamma_{\lambda}$, $P$ for $P_{\lambda}$, $G$ for $G_{\lambda}$, $\Omega$ for $\Omega_{P_{\lambda} \subseteq G_{\lambda}}$ and $O$ for $O_{P_{\lambda} \subseteq G_{\lambda}}$. Since $\Omega = O \: \dotcup \gekl{\chi_{\infty}}$,
\bgl
\label{ses-gphom}
  0 \to C_c(O,\Zz) \to C_c(\Omega,\Zz) \to \Zz \to 0
\egl
is an exact sequence of $\Zz G$-modules. It is easy to see that $\menge{1_{(x,1)P}}{x \in \Gamma}$ forms a $\Zz$-basis of $C_c(\Omega,\Zz)$ on which $\Gamma$ acts freely transitively. Hence $C_c(\Omega,\Zz) \cong \Zz \Gamma$ as $\Zz \Gamma$-modules. Thus 
$$
  H_q(\Gamma,C_c(\Omega,\Zz)) \cong
  \bfa
    \Zz & if \ q=0, \\
    \gekl{0} & if \ q \geq 1.
  \efa
$$
The Hochschild-Serre spectral sequence (see \cite[Chapter~VII, Theorem~(6.3)]{Brown}) then yields
$$
  H_i(G,C_c(\Omega,\Zz)) \cong H_i(\spkl{\lambda}) \cong
  \bfa
    \Zz & if \ i=0,1, \\
    \gekl{0} & if \ i \geq 2.
  \efa
$$
Moreover, we know that $H_0(G) \cong \Zz$ and $H_1(G) \cong H_0(\spkl{\lambda},H_1(\Gamma)) \oplus H_1(\spkl{\lambda}) \cong \Gamma / (1 - \lambda) \Gamma \oplus \spkl{\lambda} \cong \Gamma / (1 - \lambda) \Gamma \oplus \Zz$. Under these isomorphisms, $H_1(G,C_c(\Omega,\Zz)) \to H_1(G)$ is given by $\Zz \to \Gamma / (1 - \lambda) \Gamma \oplus \Zz, \, z \ma (0,z)$. Plugging this into the long exact sequence
\bglnoz
  \dotso &\to& H_1(G,C_c(O,\Zz)) \to H_1(G,C_c(\Omega,\Zz)) \to H_1(G) \\
  &\to& H_0(G,C_c(O,\Zz)) \to H_0(G,C_c(\Omega,\Zz)) \to H_0(G) \to 0
\eglnoz
attached to the short exact sequence \eqref{ses-gphom}, we arrive at the desired result.
\eproof
Note that Proposition~\ref{gphom} and Lemma~\ref{criterion_Cartan} prove Theorem~\ref{THM_Q-D}.

Combining Proposition~\ref{gphom} and Proposition~\ref{gphom'}, we obtain
\bcor
\label{HGO}
If $\lambda \in (0,\infty)$ is an algebraic integer, then
\bglnoz
  H_0(G_{\lambda}, C_c(O_{P_{\lambda} \subseteq G_{\lambda}},\Zz))
  &\cong& \Gamma_{\lambda} / (1 - \lambda) \Gamma_{\lambda} \\
  H_k(G_{\lambda}, C_c(O_{P_{\lambda} \subseteq G_{\lambda}},\Zz))
  &\cong& \Coker(\id - \Lambda^{k+1} \varphi) \oplus \Ker(\id - \Lambda^k \varphi) \ \ \ for \ all \ k \geq 1.
\eglnoz
\ecor

A comparison with Corollary~\ref{K-Q-algint} yields
\bcor
Let $\lambda \in (0,\infty)$ be an algebraic integer with $\lambda \neq 1$. Then $$
  K_i(\cQ^{\lambda})
  \cong \bigoplus_{j=0}^{\infty} 
  H_{2j+i}(G_{\lambda},C_c(O_{P_{\lambda} \subseteq G_{\lambda}},\Zz)) \ \ \ for \ i=0,1.
$$
\ecor
In other words, for the groupoid $\cG = G_{\lambda} \ltimes O_{P_{\lambda} \subseteq G_{\lambda}}$, we obtain $K_i(C^*_r(\cG)) \cong \bigoplus_{j=0}^{\infty} H_{2j+i}(\cG)$, as in \cite[\S~3.1]{Matui}.

\bremark
A similar, but much simpler computation gives for every subgroup $\Gamma$ of $\Rz$ that $H_k(\Gamma,C_c(O_{\Gamma^+ \subseteq \Gamma},\Zz)) \cong H_{k+1}(\Gamma)$, and $K_i(E(J_{\Gamma^+ \subseteq \Gamma} \rtimes \Gamma)E) \cong \bigoplus_{j=0}^{\infty} H_{2j+i}(\Gamma,C_c(O_{\Gamma^+ \subseteq \Gamma},\Zz))$ for $i=0,1$. In particular, for all $\lambda \in (0,\infty)$, $\lambda \neq 1$, we get
$$
  K_i(F^{\lambda}) 
  \cong \bigoplus_{j=0}^{\infty} 
  H_{2j+i}(\Gamma_{\lambda},C_c(O_{\Gamma_{\lambda}^+ \subseteq \Gamma_{\lambda}},\Zz))
  \ \ \ for \ i=0,1.
$$
\eremark

\subsection{Concrete computations}
\label{ConCom}

Let $1 \neq \lambda \in (0,\infty)$ be an algebraic integer with minimal polynomial $f(T) = T^d + a_{d-1} T^{d-1} + \dotso + a_1 T + a_0$. Note that irreducibility of $f$ poses conditions on $a_{d-1}, \dotsc, a_0$, for instance $a_0 \neq 0$ and $f(1) = 1 + a_{d-1} + \dotso + a_0 \neq 0$ ($\lambda \neq 1$). Because of Corollary~\ref{K-Q-algint} and Corollary~\ref{HGO}, all we have to do is to determine kernel and cokernel of $\id - \Lambda^k \varphi$. With respect to the canonical $\Zz$-basis $1, \lambda, \dotsc, \lambda^{d-1}$ of $\Zz[\lambda]$, $\varphi$ is given by the matrix
$$
\rukl{
  \begin{array}{cccc}
  0 & \dotso & 0 & -a_0 \\
  1 & \ddots & \vdots & \vdots \\
  0 & \ddots & 0 & \vdots \\
  0 & 0 & 1 & -a_{d-1}
  \end{array}
  }.
$$
It is easy to see that $\Lambda^k \varphi$ is given by the matrix of $k$-minors. So we form the difference between the identity matrix and the matrix of $k$-minors, and compute kernel and cokernel by determining the Smith normal form. All this is straightforward. In the following, we omit the details in our computations and present the final results.

\blemma
\label{comp-general}
Let $1 \neq \lambda \in (0,\infty)$ be an algebraic integer with minimal polynomial $f(T) = T^d + a_{d-1} T^{d-1} + \dotso + a_1 T + a_0$, and let $\varphi: \: \Zz[\lambda] \to \Zz[\lambda]$ be the homomorphism given by multiplication with $\lambda$.

We have $\Ker(\id - \Lambda^1 \varphi) = \gekl{0}$ and $(\Coker(\id - \Lambda^1 \varphi),1) \cong (\Zz / f(1) \Zz,1)$.

For all $d \geq 1$, $\Ker(\id - \Lambda^{d-1} \varphi) \cong \gekl{0}$, and $\Coker(\id - \Lambda^{d-1} \varphi) \cong \Zz / \tfrac{f((-1)^d a_0)}{a_0} \Zz$. Note that $a_0 \neq 0$ and that $f((-1)^d a_0) \neq 0$ since $f$ is irreducible.

Also, for all $d \geq 1$, we have
$$
\Ker(\id - \Lambda^d \varphi) \cong
\bfa
  \gekl{0} & if \ 1 + (-1)^{d+1} a_0 \neq 0, \\
  \Zz & if \ 1 + (-1)^{d+1} a_0 = 0,
\efa
\ \ \
and \ \Coker(\id - \Lambda^{d-1} \varphi) \cong \Zz / (1 + (-1)^{d+1} a_0) \Zz.$$
\elemma
Note that we write $1$ for the class of $1$ in the corresponding quotient.

We may now combine Corollary~\ref{K-Q-algint}, Corollary~\ref{HGO} and Lemma~\ref{comp-general}. Here are some concrete examples:

\begin{tabular}{c||c|c||c}
& \multicolumn{2}{ c|| }{$H_k(G_{\lambda}, C_c(O_{P_{\lambda} \subseteq G_{\lambda}},\Zz))$} & $(K_0(\cQ^{\lambda}),[1]_0,K_1(\cQ^{\lambda}))$ \\ \hline \hline
\multicolumn{1}{ c|| }{\multirow{2}{*}{\scriptsize $d=1$}} & {\scriptsize $\Zz / (1 + a_0) \Zz$} & {\scriptsize $k=0$} & {\multirow{2}{*}{\scriptsize $(\Zz / (1 + a_0) \Zz,1,\gekl{0})$}} \\ \cline{2-3}
\multicolumn{1}{ c|| }{} & {\scriptsize $\gekl{0}$} & {\scriptsize $k \geq 1$} &     \\ \hline \hline
\multicolumn{1}{ c|| }{\multirow{3}{*}{\scriptsize $d=2, a_0 = 1$}} & {\scriptsize $\Zz / (2 + a_1) \Zz$} & {\scriptsize $k=0$} & {\multirow{3}{*}{\scriptsize $(\Zz / (2 + a_1) \Zz \oplus \Zz,(1,0),\Zz)$}} \\ \cline{2-3}
\multicolumn{1}{ c|| }{} & {\scriptsize $\Zz$} & {\scriptsize $k=1,2$} &     \\ \cline{2-3}
\multicolumn{1}{ c|| }{} & {\scriptsize $\gekl{0}$} & {\scriptsize $k \geq 3$} &     \\ \hline \hline
\multicolumn{1}{ c|| }{\multirow{3}{*}{\scriptsize $d=2, a_0 \neq 1$}} & {\scriptsize $\Zz / (1 + a_1 + a_0) \Zz$} & {\scriptsize $k=0$} & {\multirow{3}{*}{\scriptsize $(\Zz / (1 + a_1 + a_0) \Zz,1,\Zz / (1 - a_0) \Zz)$}} \\ \cline{2-3}
\multicolumn{1}{ c|| }{} & {\scriptsize $\Zz / (1 - a_0) \Zz$} & {\scriptsize $k=1$} &     \\ \cline{2-3}
\multicolumn{1}{ c|| }{} & {\scriptsize $\gekl{0}$} & {\scriptsize $k \geq 2$} &     \\ \hline \hline
\multicolumn{1}{ c|| }{\multirow{4}{*}{\scriptsize $d=3, a_0 = -1$}} & {\scriptsize $\Zz / (a_2 + a_1) \Zz$} & {\scriptsize $k=0$} & {\multirow{4}{*}{\scriptsize $(\Zz / (a_2 + a_1) \Zz \oplus \Zz,(1,0),\Zz / (a_2 + a_1) \Zz \oplus \Zz)$}} \\ \cline{2-3}
\multicolumn{1}{ c|| }{} & {\scriptsize $\Zz / (a_2 + a_1) \Zz$} & {\scriptsize $k=1$} &     \\ \cline{2-3}
\multicolumn{1}{ c|| }{} & {\scriptsize $\Zz$} & {\scriptsize $k=2,3$} &     \\ \cline{2-3}
\multicolumn{1}{ c|| }{} & {\scriptsize $\gekl{0}$} & {\scriptsize $k \geq 4$} &     \\ \hline \hline
\multicolumn{1}{ c|| }{\multirow{4}{*}{\scriptsize $d=3$, $a_0 \neq -1$}} & {\scriptsize $\Zz / (1 + a_2 + a_1 + a_0) \Zz$} & {\scriptsize $k=0$} & {\multirow{4}{8cm}{\scriptsize $(\Zz / (1 + a_2 + a_1 + a_0) \Zz \oplus \Zz / (1 + a_0) \Zz,(1,0)$, $\Zz / (- a_0^2 + a_0a_2 - a_1 + 1) \Zz)$}} \\ \cline{2-3}
\multicolumn{1}{ c|| }{} & {\scriptsize $\Zz / (- a_0^2 + a_0a_2 - a_1 + 1) \Zz$} & {\scriptsize $k=1$} &     \\ \cline{2-3}
\multicolumn{1}{ c|| }{} & {\scriptsize $\Zz / (1 + a_0) \Zz$} & {\scriptsize $k=2$} &     \\ \cline{2-3}
\multicolumn{1}{ c|| }{} & {\scriptsize $\gekl{0}$} & {\scriptsize $k \geq 3$} &     \\
\end{tabular}

Note that $\lambda = -a_0$ for $d=1$, so that $a_0 \neq -1$.

Here is an explicit example: Let $\lambda = \sqrt{n}$, where $n \in \Nz$ is not a square, so that the minimal polynomial of $\lambda$ is $f(T) = T^2 - n$. Looking at the case $d=2$, $a_0 = -n \neq 1$, we obtain ($a_1 = 0$):
\bglnoz
  && H_k(G_{\sqrt{n}}, C_c(O_{P_{\sqrt{n}} \subseteq G_{\sqrt{n}}},\Zz))
  \cong
  \bfa
    \Zz / (1-n) \Zz & k=0, \\
    \Zz / (1+n) \Zz & k=1, \\
    \gekl{0} & k \geq 2;
  \efa
  \\
  && (K_0(\cQ^{\sqrt{n}}),[1]_0,K_1(\cQ^{\sqrt{n}})) \cong (\Zz / (1-n) \Zz,1,\Zz / (1+n) \Zz).
\eglnoz
Comparing this with the example in \cite[\S~2.4]{CPPR}, we see that the computation of $K_1(\cQ^{\tfrac{1}{\sqrt{n}}}) \cong K_1(\cQ^{\sqrt{n}})$ in \cite{CPPR} is wrong due to a mistake in the computation of the homomorphism $\Lambda^2 \varphi$.

Let us now comment on the examples appearing in \cite[Proposition~2.22~(2)]{CPPR}: If the minimal polynomial of $0 < \lambda < 1$ is $f(T) = T^3 + (n+1)T^2 + nT + 1$ for $n \leq -2$, then looking at our computation for $d=3$, $a_0 \neq -1$, we obtain
$$
  (K_0(\cQ^{\lambda}),[1]_0,K_1(\cQ^{\lambda})) \cong (\Zz / (2n+3) \Zz \oplus \Zz / 2 \Zz, (1,0), \gekl{0}) \cong (\Zz / (4n+6) \Zz, 2, \gekl{0}).
$$
If the minimal polynomial of $0 < \lambda < 1$ is $f(T) = T^3 + (n-1)T^2 + nT + 1$ for $n \leq -1$, then looking at our computation for $d=3$, $a_0 \neq -1$, we obtain
$$
  (K_0(\cQ^{\lambda}),[1]_0,K_1(\cQ^{\lambda})) \cong (\Zz / (2n+1) \Zz \oplus \Zz / 2 \Zz, (1,0), \gekl{0}) \cong (\Zz / (4n+2) \Zz, 2, \gekl{0}).
$$
So in either case, we see that $\cQ^{\lambda}$ is only stably isomorphic to a Cuntz algebra (of the form $\cO_{4k+3}$, $k \geq 0$), but not isomorphic as $[1]_0$ does not generate $K_0(\cQ^{\lambda})$.

Moreover, in the examples in \cite[Proposition~2.22~(2)]{CPPR}, both $\lambda$ and $\lambda^{-1}$ are algebraic integers. We now show that if we just require $\lambda$ to be an algebraic integer, then all Cuntz algebras appear. Namely, for every $n = 2,3,\dotsc$, we can find a quadratic integer $\lambda$ such that $\cQ^{\lambda} \cong \cO_n$ (as unital C*-algebras, not only up to stabilisation). To see this, take $n = 2,3,\dotsc$ and set $a \defeq -2-n$. Let $f(T) \defeq T^2 + aT + 2$, and let $\lambda \defeq - \tfrac{a}{2} + \sqrt{\tfrac{a^2}{4} - 2}$ be a root of $f$. Obviously, $\lambda > 0$, and $\lambda \in \Rz \setminus \Qz$, so that $f$ is the minimal polynomial of $\lambda$. From our computations for $d=2$, $a_0 \neq 1$, we obtain (for $a_1 = a$, $a_0 = 2$):
$$
  (K_0(\cQ^{\lambda}),[1]_0,K_1(\cQ^{\lambda})) \cong (\Zz / (1-n) \Zz, 1, \gekl{0}) \cong (K_0(\cO_n),[1]_0,K_1(\cO_n)).
$$

Now we come to our example of $\lambda$, $\mu$ such that $\cQ^{\lambda} \cong \cQ^{\mu}$ but $(\cQ^{\lambda},\cD^{\lambda}) \ncong (\cQ^{\mu},\cD^{\mu})$. Let $f(T) = T^2 - 3T + 1$. Since $f(0) = 1$, $f(1) = -1$, there exists $0 < \lambda < 1$ with $f(\lambda) = 0$. $f$ is irreducible by Gauss' Lemma since $f(1) \neq 0 \neq f(-1)$. Our computation for $d=2$, $a_0 = 1$ gives in this case ($a_1 = -3$):
\bgloz
  H_k(G_{\lambda}, C_c(O_{P_{\lambda} \subseteq G_{\lambda}},\Zz))
  \cong
  \bfa
    \gekl{0} & k=0 \ or \ k \geq 3 \\
    \Zz & k=1,2; \\
  \efa
  \ \ \ \ \
  and \ \ \ \ \ (K_0(\cQ^{\lambda}),[1]_0,K_1(\cQ^{\lambda})) \cong (\Zz,0,\Zz).
\egloz
Choose $a \in \Zz$ arbitrary and let $g(T) = T^3 + aT^2 + (1-a)T - 1$. Since $g(0) = -1$, $g(1) = 1$, there exists $0 < \mu < 1$ with $g(\mu) = 0$. $g$ is irreducible by Gauss' Lemma since $g$ does not vanish in $-1$ or $1$. Our computation for $d=3$, $a_0 = -1$ gives in this case ($a_2 = a$, $a_1 = 1-a$):
\bgloz
  H_k(G_{\mu}, C_c(O_{P_{\mu} \subseteq G_{\mu}},\Zz))
  \cong
  \bfa
    \gekl{0} & k=0,1, \ or \ k \geq 4 \\
    \Zz & k=2,3; \\
  \efa
  \ \ \ 
  and \ \ \ (K_0(\cQ^{\mu}),[1]_0,K_1(\cQ^{\mu})) \cong (\Zz,0,\Zz).
\egloz
In particular, we could take $a=1$, and $g(T) = T^3 + T^2 - 1$. So indeed, by the Kirchberg-Phillips classification theorem \cite[Chapter~8]{Ror} and Lemma~\ref{criterion_Cartan}, we have found $\lambda$, $\mu$ such that $\cQ^{\lambda} \cong \cQ^{\mu}$ but $(\cQ^{\lambda},\cD^{\lambda}) \ncong (\cQ^{\mu},\cD^{\mu})$.

Finally, we remark that for every algebraic integer $\lambda$ with $1 \neq \lambda \in (0,\infty)$, $\cQ^{\lambda} \cong \cO_n$ implies that
$$
  H_k(G_{\lambda}, C_c(O_{P_{\lambda} \subseteq G_{\lambda}},\Zz))
  \cong
  \bfa
    \Zz / (n-1) \Zz & k=0, \\
    \gekl{0} & k \geq 1.
  \efa
$$
In other words, group homology does not allow us to show that $(\cQ^{\lambda},\cD^{\lambda}) \ncong (\cO_n,\cD_n)$ (if this is the case at all). Here $\cD_n$ is the canonical Cartan subalgebra of $\cO_n$.

\section{Open questions}
\label{OQ}

Our first remark is that for the groupoid $\cG = G_{\lambda} \ltimes O_{P_{\lambda} \subseteq G_{\lambda}}$, our computation shows that $K_i(C^*_r(\cG)) \cong \bigoplus_{j=0}^{\infty} H_{2j+i}(\cG)$ for $i=0,1$, at least for algebraic integers $\lambda$. Such an isomorphism holds for other types of groupoids as well (see \cite[\S~3.1]{Matui}). So a natural question would be how general this phenomenon is, and whether there is a conceptual explanation.

Secondly, it was observed in \cite{CPPR} that for algebraic numbers $\lambda$, we must have $\rk(K_0(\cQ^{\lambda})) = \rk(K_1(\cQ^{\lambda}))$. Is this the only restriction? In other words, given a UCT Kirchberg algebra $A$ with finitely generated K-groups and $\rk(K_0(A)) = \rk(K_1(A))$, can we always find an algebraic number $\lambda$ such that $\cQ^{\lambda} \cong A$? Note that in order to have $\cQ^{\lambda} \cong A$ for an algebraic integer $\lambda$, we at least have to require that $[1]_0$ generates a direct summand in $K_0(A)$.

Finally, we mention that for transcendental $\lambda \in (0,\infty)$, we obtain $H_k(G_{\lambda}, C_c(O_{P_{\lambda} \subseteq G_{\lambda}},\Zz)) \cong \bigoplus_{i=0}^{\infty} \Zz$ for all $k = 0,1,\dotsc$, and $(K_0(\cQ^{\lambda}),[1]_0,K_1(\cQ^{\lambda})) \cong (\bigoplus_{i=0}^{\infty} \Zz,e_0,\bigoplus_{i=0}^{\infty} \Zz)$, where $e_0 = (1,0,0,\dotsc) \in \bigoplus_{i=0}^{\infty} \Zz$. This K-theory computation, regardless of the position of $[1]_0$, has been carried out in \cite[Proposition~2.24]{CPPR}. Again, we see that for transcendental $\lambda$, $\cQ^{\lambda}$ is stably isomorphic, but not isomorphic to $\cQ_{\Nz}$, since $(K_0(\cQ_{\Nz}),[1]_0,K_1(\cQ_{\Nz})) \cong (\bigoplus_{i=0}^{\infty} \Zz,0,\bigoplus_{i=0}^{\infty} \Zz)$. Having studied the case of algebraic integers and transcendental numbers, it remains to study the case of algebraic $\lambda$ which are not algebraic integers. For instance, a concrete question would be: We have shown that for an algebraic number $\lambda$ with minimal polynomial $f(T) = T^d + a_{d-1}T^{d-1} + \dotso + a_0$, where $a_{d-1}, \dotsc, a_0 \in \Zz$ (so that $\lambda$ is actually an algebraic integer), $[1]_0 \in K_0(\cQ^{\lambda})$ is a torsion element with torsion order $\abs{f(1)} = \abs{1 + a_{d-1} + \dotso + a_0}$. In general, given an algebraic number $\lambda$ with minimal polynomial $f(T) = a_dT^d + a_{d-1}T^{d-1} + \dotso + a_0$, where $a_d, \dotsc, a_0 \in \Zz$, is $[1]_0 \in K_0(\cQ^{\lambda})$ again a torsion element with torsion order $\abs{f(1)}$? Note that our methods immediately show that this is true for $d=1$ (see also \cite[Proposition~2.17]{CPPR}). One might also wonder how the generalisation of Lemma~\ref{comp-general} looks like when we go from algebraic integers to arbitrary algebraic numbers.

\end{document}